\documentclass[review,preprint]{elsarticle}

\usepackage[a4paper,total={7in,8in}]{geometry}
\usepackage{framed,multirow}
\usepackage{subfigure}
\usepackage{graphicx}%
\usepackage{multirow}%
\usepackage{amsmath,amssymb,amsfonts}%
\usepackage{amsthm}%
\usepackage{mathrsfs}%
\usepackage[title]{appendix}%
\usepackage{xcolor}%
\usepackage{textcomp}%
\usepackage{manyfoot}%
\usepackage{booktabs}%
\usepackage{algorithm}%
\usepackage{algorithmicx}%
\usepackage{algpseudocode}%
\usepackage{listings}%
\usepackage{array}
\usepackage{amssymb}
\usepackage{latexsym}
\usepackage{hyperref}
\usepackage{marvosym}

\theoremstyle{thmstyleone}%
\theoremstyle{thmstyletwo}%

\theoremstyle{thmstylethree}%
\usepackage{url}
\usepackage{xcolor}
\definecolor{newcolor}{rgb}{.8,.349,.1}

\begin{document}


\begin{frontmatter}

\title{A fully-decoupled second-order-in-time and unconditionally energy stable scheme for a phase-field model of two phase flow with variable density}%

\author[1]{Jinpeng Zhang}

\author[1]{Li Luo\corref{cor1}}
\cortext[cor1]{Corresponding author: liluo@um.edu.mo}

\author[2,3]{Xiaoping Wang\corref{cor2}}
\cortext[cor2]{Corresponding author: wangxiaoping@cuhk.edu.cn}

\address[1]{Department  of  Mathematics,  University  of  Macau,  Macao  SAR,  China}
\address[2]{School of Science and Engineering, The Chinese University of Hong Kong,
Shenzhen, Guangdong 518172, China $ \&$ Shenzhen International Center for Industrial and Applied
Mathematics, Shenzhen Research Institute of Big Data, Guangdong 518172, China}
\address[3]{Department of Mathematics, the Hong Kong University of Science and Technology, Hong Kong, China}


\begin{abstract}
In this paper, we develop a second-order, fully decoupled, and energy-stable numerical scheme for the Cahn-Hilliard-Navier-Stokes model for two phase flow with variable density and viscosity. We propose a new decoupling Constant Scalar Auxiliary Variable (D-CSAV) method which is easy to generalize to schemes with high order accuracy in time.  The method is designed using the ``zero-energy-contribution'' property while maintaining conservative time discretization for the ``non-zero-energy-contribution'' terms.  A new set of scalar auxiliary variables is introduced to develop second-order-in-time, unconditionally energy stable, and decoupling-type numerical schemes. We also introduce a stabilization parameter $\alpha$ to improve the stability of the scheme by slowing down the dynamics of the scalar auxiliary variables. Our algorithm simplifies to solving three independent linear elliptic systems per time step, two of them with constant coefficients.  The update of all scalar auxiliary variables is explicit and decoupled from solving the phase field variable and velocity field.  We rigorously prove unconditional energy stability of the scheme and perform extensive benchmark simulations to demonstrate accuracy and efficiency of the method.

\end{abstract}

\end{frontmatter}


\section{Introduction}

Phase field fluid models are increasingly important in the study of multiphase flows \cite{anderson1998diffuse, lamorgese2011phase}. They offer a remarkable ability to capture complex interface dynamics and handle topological changes efficiently. Furthermore, by achieving thermodynamic consistency and upholding the energy dissipation law, these models become particularly well-suited for a wide range of multiphase flow applications. This study focuses on developing  efficient numerical techniques for the Cahn-Hilliard-Navier-Stokes (CH-NS) model proposed by Abels et al. \cite{abels2012thermodynamically} for the two-phase incompressible flow problem with varying density and viscosity. The notable advantages of this model include satisfying thermodynamic consistency and the energy law. Furthermore, in contrast to the quasi-incompressible two-phase flow model \cite{lowengrub1998quasi}, the velocity field of the mixture in this model remains divergence-free \cite{hosseini2017isogeometric}.

Several critical considerations must be taken into account when developing efficient numerical methods for Abels' model.  Firstly, it has an inherent characteristic that energy functionals decay over time. Preserving this property at the discrete level is essential to ensure the energy stability of the numerical scheme.  Secondly, achieving high-order accuracy in time is vital for simulations to accurately capture complex dynamics, especially when the two-phase interface experiences rapid and intricate changes. Thirdly, enhancing computational efficiency requires a strategic approach. System decoupling stands out as an effective method, involving the breakdown of the complex system into simpler subproblems. This decoupling strategy reduces the complexity of the numerical algorithm and the computational resources needed, ultimately decreasing overall costs.

Over the past decade, numerous first-order energy-stable numerical methods have been developed and effectively applied in various phase-field simulations \cite{grun2014two,garcke2016stable,shen2015decoupled, zhu2019efficient,chen2022highly}. {Grün and Klingbeil \cite{grun2014two}  proposed an implicit, coupled numerical scheme and demonstrated the energy stability of the scheme. This approach requires solving a coupled nonlinear system at every time step. Garcke et al. \cite{garcke2016stable} introduced a more efficient algorithm, which is a linear, coupled, and energy-stable numerical scheme utilizing the convex splitting method. Decoupling-type numerical schemes, compared to coupled schemes, can significantly improve computational efficiency.} Shen and Yang \cite{shen2015decoupled} devised two classes of fully decoupled, unconditionally energy-stable schemes for two-phase incompressible flows by employing stabilization and convex splitting techniques. These schemes were extended to tackle the moving contact line problem \cite{gao2014, LUO2017233, yu2017numerical}. In their work \cite{shen2018}, Shen et al. introduced the scalar auxiliary variable (SAV) method as an efficient numerical approach to solve gradient flow systems (see also \cite{liu2020exponential,huang2022new,li2022fully,huang2020highly,liu2024novel,jiang2022improving}). Subsequently, Zhu \cite{zhu2019efficient, zhu2020fully} integrated the SAV method with pressure stabilization techniques to address the phase-field fluid model, showing improved accuracy through numerical experiments. However, these methods involve solving the CH-NS equations with variable coefficients at each time step because of discretizing the nonlinear convective term semi-implicitly, which can result in significant computational expenses. To address this challenge, Chen and Yang \cite{chen2022highly} proposed the decoupled auxiliary scalar variable (DSAV) method, which leads to a fully decoupled and more efficient numerical scheme.

Increasing the order of accuracy in time while upholding crucial energy stability poses a significant challenge, requiring a delicate balance of intricate considerations and the resolution of complex scheme formulations. {Many second-order accurate and energy-stable numerical schemes have recently been developed \cite{khanwale2022fully,yang2019unconditionally,gong2018fully,fu2021linear}.}
Khanwale et al. \cite{khanwale2022fully} introduced a fully coupled, implicit numerical method for Abels' model, achieving second-order temporal accuracy and conditional energy stability. However, this innovative method entails high computational complexity as it requires solving coupled non-linear systems at each time step.
Yang and Dong \cite{yang2019unconditionally} devised an unconditionally energy-stable scheme for this model, in which the coupling of velocity and pressure requires a subiteration procedure at each time step.
Fu and Han \cite{fu2021linear} developed a linear,  coupled second-order-in-time, unconditionally energy-stable numerical scheme based on the SAV method. {This method requires solving a fully coupled linear system at each time step, leading to high computational costs.} In \cite{guermond2011error},  Guermond et al. formulated a second-order-in-time and energy-stable numerical scheme for the Navier-Stokes equations with variable density. {While this discrete method can achieve unconditional energy stability, in practice, a large number of variable coefficients introduce unpredictable effects on the condition number of the linear system. It is shown in \cite{chen2022highly}, through numerical examples, that this scheme leads to increased computational costs}. Furthermore, this research only established the energy stability of the standard penalty method, which is only first-order time accuracy in pressure.

Developing an unconditionally energy-stable scheme with full decoupling and second-order time accuracy for the two-phase flow model, especially with variable density, poses significant challenges.  In the derivation of the energy law for the phase field fluid model, it is observed in \cite{chen2022highly} that many nonlinear coupling terms exhibit a ``zero-energy-contribution (ZEC)'' feature. A decoupling method was introduced in \cite{chen2022highly} to address these ZEC terms. However, extending this approach to a second-order numerical scheme presents difficulties, particularly in discretizing terms related to the time derivative in the Navier-Stokes equations. { Wang et al. \cite{wang2024efficient} introduced a decoupled, linear, second-order in time, and energy-stable numerical scheme. The artificial compressibility scheme was employed to decouple pressure from velocity. Nonetheless, this approach has the drawback of introducing non-physical acoustic waves, resulting in oscillations in the presence of a large density ratio.}
 
It is important to note that in addition to the ZEC terms, there are also ``non-zero-energy-contribution (NZEC)'' terms which must be time conservative, meaning that they can be expressed as time derivatives of an energy functional (e.g. the kinetic energy) in the derivation of the energy law. However, this conservative structure can be compromised following the time discretization of the PDE system. In this paper, we present a general procedure to address these nonzero energy-contribution terms.
In a recent investigation of gradient flow problems \cite{zhang2024new}, a constant scalar auxiliary variable and its corresponding ODE were introduced. This established a connection between the conservative form and the non-conservative form of the nonlinear free energy after time discretization. Through this auxiliary variable, the non-conservative form can be substituted by the conservative form in the energy proof, ensuring unconditional stability for the discrete energy. Drawing inspiration from \cite{zhang2024new}, we introduce a  new set of constant scalar auxiliary variables and their governing ODEs to tackle the challenges associated with discretizing the time derivative terms in the Navier-Stokes equations. This strategy facilitates the development of an unconditionally energy-stable scheme with a full decoupling structure and second-order time accuracy, termed D-CSAV. Furthermore, we devise a modified rotational penalty method based on D-CSAV and establish its energy stability.  We also introduce a stabilization parameter $\alpha$ to improve the stability of the scheme by slowing down the dynamics of the scalar auxiliary variables. Extensive numerical experiments demonstrate that this method enhances the time accuracy of the pressure compared to the standard version, with convergence verified through extensive numerical testing.

{
The structure of the remainder of the article is as follows. Section 2 introduces the model and outlines the time-stepping method and rigorously proves its unconditional energy stability; Section 3 presents numerical experiments to verify the accuracy and stability of the proposed scheme; and Section 4 concludes with a brief summary and additional discussions.
}

\section{The model and discretization}
\subsection{The model}
\noindent 
We consider the Cahn-Hilliard Navier-Stokes model for a two-phase fluid flow system with variable density and viscosity as proposed by Abels et al. \cite{abels2012thermodynamically}. This model adheres to thermodynamic consistency and the energy law, ensuring the velocity field remains solenoidal \cite{hosseini2017isogeometric}. The domain $\Omega \subset \mathbf{R}^{d}$ ($d = 2,3$) containing a two-phase immiscible fluid mixture is a smooth, rectangular, open, bounded, and connected region.
The governing equations are given by
\begin{align}
& \phi_t + \nabla \cdot (\boldsymbol{u} \phi) = \nabla \cdot \left(M(\phi) \nabla \mu \right), \label{CHNS_1} \\
& \mu = \lambda \left(-\epsilon \Delta \phi + \frac{1}{\epsilon} f(\phi) \right), \label{CHNS_2} \\
& \rho(\phi) \left(\boldsymbol{u}_t + (\boldsymbol{u} \cdot \nabla) \boldsymbol{u} \right) + \boldsymbol{J}(\mu) \cdot \nabla \boldsymbol{u} - \nabla \cdot (\nu(\phi) D(\boldsymbol{u})) + \nabla p + \phi \nabla \mu = \rho \mathbf{g}, \label{CHNS_3} \\
& \nabla \cdot \boldsymbol{u} = 0, \label{CHNS_4}
\end{align}
with the following definitions:
\begin{align}
& \boldsymbol{J}(\mu) = M(\phi) \nabla \mu \frac{\rho_2 - \rho_1}{2}, \label{CHNS_5} \\
& D(u) = \nabla \boldsymbol{u}^T + \nabla \boldsymbol{u}, \label{CHNS_6} \\
& \rho(\phi) = \frac{\phi}{2} (\rho_1 - \rho_2) + \frac{1}{2} (\rho_1 + \rho_2), \label{CHNS_7} \\
& \nu(\phi) = \frac{\phi}{2} (\nu_1 - \nu_2) + \frac{1}{2} (\nu_1 + \nu_2). \label{CHNS_8}
\end{align}
Here, $\phi$ represents the phase-field variable that describes the transition between the two homogeneous equilibrium phases $\phi_{\pm} = \pm 1$. $\mu$ is the chemical potential, $\boldsymbol{u}$ is the fluid velocity field, $p$ is the pressure, $\rho(\phi)$ is the density, $\nu(\phi)$ is the viscosity, and $\rho_{j}$ and $\nu_{j} > 0$ are the specific (constant) density and viscosity of fluid $j = 1, 2$. The function $f(\phi) = \phi^3 - \phi$ represents the derivative of the nonlinear double-well potential $F(\phi) = \frac{1}{4}(\phi^2 - 1)^2$.
Moreover, $D(\boldsymbol{u}) = \left(\nabla \boldsymbol{u} + \nabla \boldsymbol{u}^T\right)$ denotes the fluid strain tensor, $\boldsymbol{J}(\mu)$ represents the relative flux related to component diffusion, $M(\phi)$ is the mobility coefficient (which can be a positive constant or a positive function related to $\phi$), $\lambda$ is the surface tension parameter, $\epsilon$ is a small parameter related to the interfacial region's thickness, and $\mathbf{g} =  g_{0}\mathbf{\hat{g}}$ where $g_{0}$ is the gravitational constant, and $\mathbf{\hat{g}}$ is a unit vector denoting the direction of gravity.

The initial conditions and the boundary conditions are given by
\begin{align}
& \boldsymbol{u}|_{t=0} = \boldsymbol{u}^{0}(x), \quad p|_{t=0} = p^{0} = 0, \quad \phi|_{t=0} = \phi^{0}(x), \label{CHNS_init}\\
& \boldsymbol{u}|_{\partial\Omega} = \mathbf{0}, \quad \partial_{\boldsymbol{n}}\phi|_{\partial\Omega} = \partial_{\boldsymbol{n}}\mu|_{\partial\Omega} = 0, \label{CHNS_bounda}
\end{align}
where $\mathbf{n}$ is the unit outward normal on the boundary $\partial\Omega$. Periodic boundary conditions can also be applied in this model.

The following theorem gives the energy dissipation law of the system when the external force is ignored ($\mathbf{g}=\mathbf{0}$). We define the following notations: $(f,g)$ means the $L^{2}$ inner product of any two functions $f(\boldsymbol{x})$ and $g(\boldsymbol{x})$, $\int_{\Omega}f(\boldsymbol{x})g(\boldsymbol{x})d\boldsymbol{x}$. $\|f(\boldsymbol{x})\|$ means the $L^{2}$ norm of $f(\boldsymbol{x})$. \\
$\textbf{Theorem 1.}$ Ignoring the gravitational term $\rho\mathbf{g}$ in \eqref{CHNS_3}, the solutions of the system \eqref{CHNS_1}-\eqref{CHNS_8} with the initial conditions \eqref{CHNS_init} and  
 the boundary conditions \eqref{CHNS_bounda} satisfy the following law of energy dissipation,
\begin{align}
\frac{d}{d t} E(\rho, \boldsymbol{u}, \phi)=-\int_{\Omega}M(\phi)|\nabla \mu|^2 d \boldsymbol{x}-\frac{1}{2} \int_{\Omega} \nu|D(\boldsymbol{u})|^2 d \boldsymbol{x} \leqslant 0, \label{EO_law}
\end{align}
where
\begin{align}
E(\rho, \boldsymbol{u}, \phi)=\int_{\Omega}\left(\frac{1}{2} \rho|\boldsymbol{u}|^2+\lambda \frac{\epsilon}{2}|\nabla \phi|^2+\frac{\lambda}{\epsilon} F(\phi)\right) d \boldsymbol{x} \label{EO}
\end{align}
is the total energy of the system.\\
\textbf{Proof.} Multiplying \eqref{CHNS_1} by $\mu$ and performing integration by parts, we have
\begin{align}
\left(\phi_t, \mu\right)+\int_{\Omega}M(\phi)|\nabla \mu|^2 d \boldsymbol{x}=-\int_{\Omega} \nabla \cdot(\boldsymbol{u} \phi) \mu d \boldsymbol{x}. \label{CHNS_1_v1}
\end{align}
By taking the $L^{2}$ inner product of \eqref{CHNS_2} with $-\phi_t$ and performing integration by parts,
\begin{align}
-\left(\mu, \phi_t\right)+\lambda \frac{d}{d t} \int_{\Omega}\left(\frac{\epsilon}{2}|\nabla \phi|^2+\frac{1}{\epsilon} F(\phi)\right) d x=0. \label{CHNS_2_v1}
\end{align}
Then, we take the inner product of \eqref{CHNS_3} with $\boldsymbol{u}$, using integration by parts and applying the divergence-free condition \eqref{CHNS_4}, we obtain
\begin{align}
	\left(\rho \boldsymbol{u}_t, \boldsymbol{u}\right) & +\frac{1}{2} \int_{\Omega} \nu|D(\boldsymbol{u})|^2 d \boldsymbol{x}+\int_{\Omega}(\phi \nabla \mu \cdot \boldsymbol{u}) d \boldsymbol{x}+\int_{\Omega} \rho(\boldsymbol{u} \cdot \nabla) \boldsymbol{u} \cdot \boldsymbol{u} d \boldsymbol{x}+\int_{\Omega} \boldsymbol{J} \cdot \nabla \boldsymbol{u} \cdot \boldsymbol{u} d \boldsymbol{x}=0.  \label{CHNS_3_v1}
\end{align}
By combining \eqref{CHNS_1}, \eqref{CHNS_5} and \eqref{CHNS_7}, we derive
\begin{align}
\rho_t+\nabla \cdot(\rho \boldsymbol{u})+\nabla \cdot \boldsymbol{J}=0. \label{CHNS_1_v2}
\end{align}
By multiplying $\frac{1}{2} \boldsymbol{u}$ to \eqref{CHNS_1_v2}, we derive
\begin{align}
\frac{1}{2} \rho_t \boldsymbol{u}+\frac{1}{2} \nabla \cdot(\rho \boldsymbol{u}) \boldsymbol{u}+\frac{1}{2} \nabla \cdot \boldsymbol{J} \boldsymbol{u}=0. \label{CHNS_1_v3}
\end{align}
Then, we take the inner product of \eqref{CHNS_1_v3} with $\boldsymbol{u}$ to derive
\begin{align}
\left(\rho_t, \frac{1}{2}|\boldsymbol{u}|^2\right)+\frac{1}{2} \int_{\Omega} \nabla \cdot(\rho \boldsymbol{u}) \boldsymbol{u} \cdot \boldsymbol{u} d \boldsymbol{x}+\frac{1}{2} \int_{\Omega} \nabla \cdot \boldsymbol{J} \boldsymbol{u} \cdot \boldsymbol{u} d \boldsymbol{x}=0. \label{CHNS_1_v4}
\end{align}
By combining \eqref{CHNS_1_v1}-\eqref{CHNS_3_v1} and \eqref{CHNS_1_v4}, we get the law of energy dissipation as follows,
$$
\frac{d}{d t}\int_{\Omega}\left(\frac{1}{2} \rho|\boldsymbol{u}|^2+\lambda \frac{\epsilon}{2}|\nabla \phi|^2+\frac{\lambda}{\epsilon} F(\phi)\right) d \boldsymbol{x}=-\int_{\Omega}M(\phi)|\nabla \mu|^2 d \boldsymbol{x}-\frac{1}{2} \int_{\Omega} \nu|D(\boldsymbol{u})|^2 d x \leqslant 0.
$$
Based on the above theorem, the original system satisfies the energy law at the continuous level. In the next section, we will develop a fully decoupled, energy-stable, and more efficient second-order-in-time numerical scheme for solving the model \eqref{CHNS_1}-\eqref{CHNS_8}.\par
\subsection{Numerical schemes}
\subsubsection{The modified system}
In \cite{chen2022highly}, the authors proposed a fully decoupled and first-order-in-time numerical scheme to solve the model \eqref{CHNS_1}-\eqref{CHNS_8}. This method introduces scalar auxiliary variables along with their respective ordinary differential equations. For ease of reference, we refer to this approach as the DSAV method in the subsequent discussion. To derive the PDE energy law \eqref{EO_law}, it is noted that certain nonlinear terms cancel out as a result of integration by parts and boundary conditions \eqref{CHNS_bounda}, as follows,
\begin{align}
	& \int_{\Omega} \nabla \cdot(\boldsymbol{u} \phi) \mu d \boldsymbol{x}+\int_{\Omega}(\phi \nabla \mu \cdot \boldsymbol{u})=0, \label{ZEC_1}\\
	& \int_{\Omega} \rho(\boldsymbol{u} \cdot \nabla) \boldsymbol{u} \cdot \boldsymbol{u} d \boldsymbol{x}+\frac{1}{2} \int_{\Omega} \nabla \cdot(\rho \boldsymbol{u}) \boldsymbol{u} \cdot \boldsymbol{u} d \boldsymbol{x}=0, \label{ZEC_2_1} \\
	& \int_{\Omega} \boldsymbol{J} \cdot \nabla \boldsymbol{u} \cdot \boldsymbol{u} d \boldsymbol{x}+\frac{1}{2} \int_{\Omega} \nabla \cdot \boldsymbol{J} \boldsymbol{u} \cdot \boldsymbol{u} d \boldsymbol{x}=0. \label{ZEC_2_2}
 \end{align}
The  equations \eqref{ZEC_1}-\eqref{ZEC_2_2} show that the contribution of all these nonlinear terms to the total free energy of the system is zero, a property termed the ``zero-energy-contribution'' in \cite{chen2022highly}. This feature is exploited to introduce the following scalar auxiliary variables $Q$, $R$, and their associated ordinary differential equations (ODEs)
\begin{align}
&\left\{\begin{array}{l}
	Q_t=\int_{\Omega}(\nabla \cdot(\boldsymbol{u} \phi) \mu+(\phi \nabla \mu) \cdot \boldsymbol{u}) d \boldsymbol{x}, \\
	\left.Q\right|_{(t=0)}=1.
\end{array}\right. \label{yang_Q}\\
&\left\{\begin{array}{l}
R_t=\int_{\Omega}\left(\rho(\boldsymbol{u} \cdot \nabla) \boldsymbol{u} \cdot \boldsymbol{u}+\frac{1}{2} \nabla \cdot(\rho \boldsymbol{u}) \boldsymbol{u} \cdot \boldsymbol{u}+(\boldsymbol{J} \cdot \nabla \boldsymbol{u}) \cdot \boldsymbol{u}+\frac{1}{2} \nabla \cdot \boldsymbol{J u} \cdot \boldsymbol{u}\right) d \boldsymbol{x}, \\
	\left.R\right|_{(t=0)}=1.
\end{array}\right. \label{yang_R}
\end{align}
It is easy to see that $Q_t=0$, $R_t=0$, so the value of the auxiliary variables $Q$ and $R$ is 1. However, in their approach, the auxiliary variables ($Q$, $R$) and the unknowns ($\boldsymbol{u}$, $\phi$, $\mu$, $p$) in the model are coupled together after numerical discretization. To establish the decoupling-type numerical schemes, the unknowns ($\boldsymbol{u}$, $\phi$, $\mu$, $p$) need to be split into a linear combination form consisting of the nonlocal auxiliary variables ($Q$, $R$). Therefore, the system \eqref{CHNS_1}-\eqref{CHNS_4} is also split into seven independent sub-systems during numerical calculations, leading to relatively high computational cost.\par
{ Notice  that in the proof of  ${Theorem 1}$, there are nonlinear terms such as $\rho \boldsymbol{u}_t$, $\rho_t \boldsymbol{u}$, and $f(\phi)$ that contribute to the kinetic energy and free energy. We refer them to the ``non-zero energy contribution (NZEC)'' terms. At the continuous level, these terms can be written as time derivatives of an energy functional (e.g., the free energy and the kinetic energy) in the derivation of the energy law, indicating they are time-conservative. However, this conservative structure may be compromised after the time discretization of the PDE system, posing challenges in designing an energy-stable numerical scheme. In the following,  we will discuss how to address these ``non-zero energy contribution'' terms. A new set of scalar auxiliary variables is introduced to develop second-order-in-time, unconditionally energy stable and decoupling-type numerical schemes.} We also introduce a stabilization parameter $\alpha$ to improve the stability of the scheme by slowing down the dynamics of the scalar auxiliary variables.\par

In \cite{zhang2024new}, we introduced the following constant scalar auxiliary variable (CSAV) method to solve the gradient flow system
\begin{align}
	& \phi_t=\Delta\mu, \label{Cphi}\\
	& \mu=-\Delta\phi+s\phi+r(f(\phi)-s\phi), \label{Cmu}\\
	& \frac{dr}{dt}=\alpha\big(-\frac{dE_{0}(\phi)}{dt}+r\int_{\Omega}(f(\phi)-s\phi)\phi_{t}d\Omega\big). \label{CSAV_r}
\end{align}
where $E_{0}(\phi)=\int_{\Omega}(F(\phi)-\frac{s}{2}\phi^{2})d\Omega$, $s$ is a stabilization parameter, and $r|_{(t=0)}=1$. $\alpha$ is the stabilization parameter and is usually taken to be small, which significantly slows down the change of $r$ so that its value is close to $1$, leading to better consistency and stability of the scheme. A first-order energy stable scheme is then given by
\begin{align}
	& \frac{\phi^{n+1}-\phi^{n}}{\delta t}=-\Delta \mu^{n+1}, \label{NuSCH1}\\
	& \mu^{n+1}=\delta \phi^{n+1}+s\phi^{n+1}+r^{n}(f(\phi^{n})-s\phi^{n}), \label{NuSCH2}\\
	& \frac{r^{n+1}-r^{n}}{\delta t}=\alpha\bigg(-\frac{E_{0}(\phi^{n+1})-E_{0}(\phi^{n})}{\delta t}+r^{n}\int_{\Omega}(f(\phi^{n})-s\phi^{n})\frac{\phi^{n+1}-\phi^{n}}{\delta t}d\Omega\bigg), \label{NuSCH3}
\end{align}
where $\delta t$ is the time step size, the superscript $n$ indicates the approximation of the variable at the $n^{th}$ time step. 
{\eqref{NuSCH3} establishes a connection between the conservative time discretization and the non-conservative time discretization of the nonlinear free energy $E_{0}$.  Multiplying  \eqref{NuSCH2} by $ \frac{\phi^{n+1}-\phi^{n}}{\delta t}$ and using \eqref{NuSCH3}, the third term on the right-hand side of \eqref{NuSCH2} can be replaced by the conservative form of the numerical time derivative of $E_{0}$, $\frac{E_{0}(\phi^{n+1})-E_{0}(\phi^{n})}{\delta t}$. Therefore, we can then obtain energy decay for the discrete energy.} The CSAV method also allows us to treat the nonlinear term in \eqref{Cmu} totally explicit, thereby decoupling the solution of equation \eqref{CSAV_r} from equations \eqref{Cphi} and \eqref{Cmu}. Additionally, a second-order energy-stable scheme can be derived based on CSAV, as demonstrated in \cite{zhang2024new}.

In this study, we extend the method to address the Cahn-Hilliard-Navier-Stokes equation. Five constant scalar auxiliary variables—$r, Q, R, T, K$—are introduced along with their corresponding ODEs. Consequently, the CH-NS system \eqref{CHNS_1}-\eqref{CHNS_4} is reformulated into the following equivalent system:
\begin{align}
	& \phi_t+Q \nabla \cdot(\boldsymbol{u} \phi)=\nabla\cdot \left(M(\phi)\nabla \mu\right), \label{MCHNS_1_sed}\\
	& \mu=\lambda\left(-\epsilon \Delta \phi+\frac{s}{\epsilon} \phi+\frac{r}{\epsilon}(f(\phi)-s\phi)\right), \label{MCHNS_2_sed}\\
	&K(\rho \boldsymbol{u}_t+\frac{1}{2} \rho_t \boldsymbol{u})-\nabla \cdot(\nu D(\boldsymbol{u}))+R\nabla p+Q \phi \nabla \mu \nonumber\\
	& \quad+R \rho(\boldsymbol{u} \cdot \nabla) \boldsymbol{u}+\frac{1}{2} R \nabla \cdot(\rho \boldsymbol{u}) \boldsymbol{u}+R \boldsymbol{J} \cdot \nabla \boldsymbol{u}+\frac{1}{2} R \nabla \cdot \boldsymbol{J} \boldsymbol{u}=\rho\mathbf{g}, \label{MCHNS_3_sed} \\
	& \nabla \cdot \boldsymbol{u}=0, \label{MCHNS_4_sed}\\
 & r_t=\alpha\big(-\frac{dE_{0}}{dt}+\frac{r\lambda}{\epsilon}\int_{\Omega}(f(\phi)-s\phi)\phi_{t} d \boldsymbol{x}\big), \label{MCHNS_5_sed}\\
	& Q_t=\alpha Q\int_{\Omega}(\nabla \cdot(\boldsymbol{u} \phi) \mu+(\phi \nabla \mu) \cdot \boldsymbol{u}) d \boldsymbol{x}, \label{MCHNS_6_sed}\\
	& R_t=\alpha R\int_{\Omega}\left(\rho(\boldsymbol{u} \cdot \nabla) \boldsymbol{u} \cdot \boldsymbol{u}+\frac{1}{2} \nabla \cdot(\rho \boldsymbol{u}) \boldsymbol{u} \cdot \boldsymbol{u}+(\boldsymbol{J} \cdot \nabla \boldsymbol{u}) \cdot \boldsymbol{u}+\frac{1}{2} \nabla \cdot \boldsymbol{J} \boldsymbol{u} \cdot \boldsymbol{u}+\nabla p\cdot\boldsymbol{u}\right) d \boldsymbol{x}, \label{MCHNS_7_sed}\\
	& K_t=\alpha \int_{\Omega}\left(-\frac{d}{dt}\left(\frac{1}{2}\rho|\boldsymbol{u}|^{2}\right)+K\left(\frac{1}{2}\frac{d\rho}{dt}\boldsymbol{u}+\rho\frac{d\boldsymbol{u}}{dt}\right)\cdot\boldsymbol{u}\right)d \boldsymbol{x}, \label{MCHNS_8_sed}\\
	& T_t=\alpha T\left(\int_{\Omega} p\nabla \cdot \boldsymbol{u} + \beta\nabla(\nu\nabla\cdot\boldsymbol{u})\cdot\nabla p d \boldsymbol{x}\right), \label{MCHNS_9_sed}
\end{align}
where  $E_{0}=\frac{\lambda}{\epsilon}\int_{\Omega}\left(F(\phi)-\frac{s}{2}\phi^{2}\right)d\boldsymbol{x}$. $\alpha$ is small parameter and $\beta$ is a parameter to be determined later. 
The initial conditions and the boundary conditions are given below,
$$
\left\{\begin{array}{l}
\left.\boldsymbol{u}\right|_{(t=0)}=\boldsymbol{u}^0,\left.p\right|_{(t=0)}=p^0,\left.\phi\right|_{(t=0)}=\phi^0,\left.\mu\right|_{(t=0)}=\mu^0=\lambda\left(-\epsilon \Delta \phi^0+\frac{1}{\epsilon} f\left(\phi^0\right)\right), \\
\left.r\right|_{(t=0)}=1,\left.Q\right|_{(t=0)}=1,\left.R\right|_{(t=0)}=1,\left.T\right|_{(t=0)}=1,\left.K\right|_{(t=0)}=1,\\
\boldsymbol{u}|_{\partial\Omega}=\mathbf{0}, \partial_{\boldsymbol{n}}\phi|_{\partial\Omega}=\partial_{\boldsymbol{n}}\mu|_{\partial\Omega}=0. 
\end{array}\right.
$$

We now elaborate on the auxiliary variables $r, Q, R, T, K$ and the associated ODEs in \eqref{MCHNS_1_sed}-\eqref{MCHNS_9_sed}. Firstly, it is observed that all integral terms in \eqref{MCHNS_5_sed}-\eqref{MCHNS_9_sed} amount to zero, indicating that $r(t) \equiv 1$, $Q(t) \equiv 1$, $R(t) \equiv 1$, $T(t) \equiv 1$, and $K(t) \equiv 1$. Secondly, we introduce an additional term $K\left(\frac{1}{2} \rho_t \boldsymbol{u}\right)+R\left(\frac{1}{2} \nabla \cdot(\rho \boldsymbol{u}) \boldsymbol{u}+\frac{1}{2} \nabla \cdot \boldsymbol{J} \boldsymbol{u}\right)$, which equals to zero as per \eqref{CHNS_1_v3}, into the momentum equation \eqref{MCHNS_3_sed}. 
Notice that the 
ODEs  \eqref{MCHNS_6_sed} and \eqref{MCHNS_7_sed} for $Q, R$ differ from \eqref{yang_Q} and \eqref{yang_R} (introduced  in \cite{chen2022highly}) in extra factors $\alpha Q$ and $\alpha R$ on the right hand sides. These ODEs enable treating the nonlinear term explicitly and decoupling \eqref{MCHNS_6_sed} and \eqref{MCHNS_7_sed} from the other equations (see Sections \ref{firstorder} and \ref{second-order}). 

In \cite{fu2021linear}, the authors note that the DSAV method performs reasonably well for the constant mobility scenario. However, in cases where the mobility coefficient degenerates, the method falters as the auxiliary variable $Q$ rapidly approaches zero, irrespective of the chosen small time step size. In our approach, the parameter $\alpha$ serves to boost algorithm stability and retains the auxiliary variable $Q$ in proximity to 1 (for more discussions on the effect of $\alpha$, see \cite{zhang2024new}). Notably, when faced with a degenerate mobility coefficient, as exemplified in the numerical experiments in Section 2.3.2, our method excels.

 
To deal with  the terms $\rho_t\boldsymbol{u}, \rho\boldsymbol{u}_t$ and $\nabla p$ in \eqref{MCHNS_3_sed}, two more auxiliary variables $K$ and $T$ are introduced. To derive the evolution ODEs for $K$ and $T$, notice that we have 
\begin{align}
	\int_{\Omega} \left(\rho \boldsymbol{u}_t+\frac{1}{2} \rho_t \boldsymbol{u}\right)\cdot \boldsymbol{u} d \boldsymbol{x} = \frac{d}{dt} \int_{\Omega} \left(\frac{1}{2}\rho|\boldsymbol{u}|^{2}\right) d \boldsymbol{x},
\end{align}
but this equation does not hold after the discretization in time (taking the first-order-in-time scheme for example), 
\begin{align} \label{eq40}
	\int_{\Omega}\left(\rho^{n+1} \frac{\boldsymbol{u}^{n+1}-\boldsymbol{u}^{n}}{\delta t}+\frac{1}{2} \frac{\rho^{n+1}-\rho^n}{\delta t} \boldsymbol{u}^{n+1}\right)\cdot \boldsymbol{u^{n+1}}d \boldsymbol{x}
	\neq   \frac{1}{\delta t}\int_{\Omega}{\frac{1}{2}\rho^{n+1}|\boldsymbol{u}^{n+1}|^{2}-\frac{1}{2}\rho^{n}|\boldsymbol{u}^{n}|^{2}}d \boldsymbol{x}.
\end{align}
Therefore the conservative
structure is lost after time discretization. To overcome this difficulty, we introduce the  auxiliary variables $K$ in \eqref{MCHNS_3_sed} that satisfies the trivial ODE ($K\equiv 1$) 
\begin{align}
\left\{\begin{array}{l}
K_t=\alpha \int_{\Omega}\left(-\frac{d}{dt}\left(\frac{1}{2}\rho|\boldsymbol{u}|^{2}\right)+K\left(\frac{1}{2}\frac{d\rho}{dt}\boldsymbol{u}+\rho\frac{d\boldsymbol{u}}{dt}\right)\cdot\boldsymbol{u}\right)d \boldsymbol{x}, \label{Kt}\\
	\left.K\right|_{(t=0)}=1, 
\end{array}\right. 
\end{align}
%
where $\alpha$ is a small constant to maintain the numerical stability. { Similar to Equation \eqref{NuSCH3}, Equation \eqref{Kt} establishes, following time discretization, a connection between the conservative time discretization and the non-conservative time discretization of the kinetic energy (refer to Sections \ref{firstorder} and \ref{second-order}). This connection allows for replacing the non-conservative time discretization (the left-hand side of \eqref{eq40}) with its conservative approximation (the right-hand side of  \eqref{eq40}), thereby leading to the derivation of the energy-decaying property.}

Next, we cope with the terms associated with the pressure. In \cite{zhang2024efficient}, a modified penalty method is introduced to establish an unconditionally energy-stable numerical scheme for the CH-NS model with non-homogenous boundary conditions. This method involves the introduction of the following scalar auxiliary variable $T$ to handle the difficulties caused by non-homogeneous boundary conditions when addressing the pressure term,
\begin{align}
\left\{\begin{array}{l}
	T_t=\alpha\int_{\Omega}(p \nabla \cdot u)d \boldsymbol{x}, \\
	\left.T\right|_{(t=0)}=1. 
\end{array}\right. \label{AV_T}
\end{align}
In \eqref{MCHNS_9_sed},  a modified ODE for $T$ is introduced to improve the accuracy of the pressure, which can also be expressed as
$$
\left\{\begin{array}{l}
	T_t=\alpha T\left(\int_{\Omega}(p\nabla \cdot u) + \beta\nabla(\nu\nabla\cdot\boldsymbol{u})\cdot\nabla pd \boldsymbol{x}\right), \\
	\left.T\right|_{(t=0)}=1.
\end{array}\right.
$$
More details for this ODE will be explained in the Sections \ref{firstorder} and \ref{second-order}.

\noindent $\textbf{Theorem 4.}$ Ignoring the gravitational term $\rho\mathbf{g}$ in \eqref{MCHNS_3_sed}, the solutions of the system \eqref{MCHNS_1_sed}-\eqref{MCHNS_9_sed} with the boundary condition \eqref{CHNS_bounda} and the initial condition \eqref{CHNS_init} satisfy the following law of energy dissipation,
\begin{align}
	\frac{d}{d t}E_{M}
	=-\int_{\Omega}M(\phi)|\nabla \mu|^2 d \boldsymbol{x}-\frac{1}{2} \int_{\Omega} \nu|D(\boldsymbol{u})|^2 d x \leq 0, \label{EM_decay_sed}
\end{align}
where
\begin{align}
	E_{M}=\int_{\Omega}\left(\frac{1}{2} \rho|\boldsymbol{u}|^2+\lambda \frac{\epsilon}{2}|\nabla \phi|^2+\frac{\lambda}{\epsilon}F(\phi)\right) d\boldsymbol{x}+\frac{1}{\alpha}r+\frac{1}{\alpha}Q+\frac{1}{\alpha}R+\frac{1}{\alpha}T+\frac{1}{\alpha}K. \label{EM}
\end{align}
\noindent \textbf{Proof.} By taking the inner product of \eqref{MCHNS_1_sed} with $\mu$ in $L^2$, we derive
\begin{align}
	\left(\phi_t, \mu\right)+\int_{\Omega}M(\phi)|\nabla \mu|^2 d \boldsymbol{x}=-Q \int_{\Omega} \nabla \cdot(\boldsymbol{u} \phi) \mu d\boldsymbol{x}. \label{MCHNS_1_v1_sed}
\end{align}
Then, we take the $L^2$ inner product of \eqref{MCHNS_2_sed} with $-\phi_t$, 
\begin{align}
	-\left(\mu, \phi_t\right)+\lambda \frac{d}{d t}\left(\int_{\Omega} \frac{\epsilon}{2}|\nabla \phi|^2+\frac{s}{2 \epsilon}|\phi|^2\right) d\boldsymbol{x}=-\frac{\lambda r}{\epsilon}\int_{\Omega}\left(f(\phi)-s\phi\right)\phi_{t}d\boldsymbol{x}. \label{MCHNS_2_v1_sed}
\end{align}
By taking the $L^2$ inner product of \eqref{MCHNS_3_sed} with $\boldsymbol{u}$ and using integration by parts and the divergence free condition \eqref{MCHNS_4_sed}, we obtain
\begin{align}
	K\frac{d}{d t} \int_{\Omega} \frac{1}{2} \rho|\boldsymbol{u}|^2 d \boldsymbol{x}= & -\frac{1}{2} \int_{\Omega} \nu|D(\boldsymbol{u})|^2 d \boldsymbol{x}-Q \int_{\Omega}(\phi \nabla \mu) \cdot \boldsymbol{u} d \boldsymbol{x} \label{MCHNS_3_v1_sed}\\
	& -R \int_{\Omega} \rho(\boldsymbol{u} \cdot \nabla) \boldsymbol{u} \cdot \boldsymbol{u} d \boldsymbol{x}-\frac{1}{2} R \int_{\Omega} \nabla \cdot(\rho \boldsymbol{u}) \boldsymbol{u} \cdot \boldsymbol{u} d \boldsymbol{x} \nonumber\\
	& -R \int_{\Omega}(\boldsymbol{J} \cdot \nabla \boldsymbol{u}) \cdot \boldsymbol{u} d \boldsymbol{x}-\frac{1}{2} R \int_{\Omega} \nabla \cdot \boldsymbol{J}\boldsymbol{u} \cdot \boldsymbol{u} d \boldsymbol{x}-R\int_{\Omega}\nabla p \cdot \boldsymbol{u} d \boldsymbol{x}. \nonumber
\end{align}
By combining the above equations \eqref{MCHNS_5_sed}-\eqref{MCHNS_9_sed} and \eqref{MCHNS_1_v1_sed}-\eqref{MCHNS_3_v1_sed}, we obtain the law of energy dissipation of the modified system \eqref{MCHNS_1_sed}-\eqref{MCHNS_7_sed} as follows,
\begin{align}
	&\frac{d}{d t}\left(\int_{\Omega}\left(\frac{1}{2} \rho|\boldsymbol{u}|^2+\lambda \frac{\epsilon}{2}|\nabla \phi|^2+\frac{\lambda}{\epsilon}F(\phi)\right) d\boldsymbol{x}+\frac{1}{\alpha}r+\frac{1}{\alpha}Q+\frac{1}{\alpha}R+\frac{1}{\alpha}T+\frac{1}{\alpha}K \right) \nonumber\\
	& =-\int_{\Omega}M(\phi)|\nabla \mu|^2 d \boldsymbol{x}-\frac{1}{2} \int_{\Omega} \nu|D(\boldsymbol{u})|^2 d x \leq 0. \label{MCHNS_EM_decay_sed}
\end{align}
$\textbf{Remark.}$ It is worth noting that the modified energy $E_{M}$ in \eqref{EM} and the original energy $E$ in \eqref{EO} are equivalent, differing only by a constant.
\begin{align}
	E_{M} - E= \frac{1}{\alpha}r+\frac{1}{\alpha}Q+\frac{1}{\alpha}R+\frac{1}{\alpha}T+\frac{1}{\alpha}K. \label{EOEM}
\end{align}
\subsubsection{Numerical schemes}
In this subsection, we develop both the first-order-in-time and the second-order-in-time marching  scheme to solve the system \eqref{MCHNS_1_sed}-\eqref{MCHNS_9_sed}. We consider the time domain $[0,T]$ and discretize it into equally spaced intervals $0=t_{0}<t_{1}<\cdots<t_{N}=T$, where $t_{i}=i\delta t$ and $\delta t = T/N$. Let $\psi^{n}$ represent the numerical approximation to the function $\psi(\cdot,t^n)$.
\subsubsection{The first-order scheme} \label{firstorder}
Given $\boldsymbol{u}^n$, $p^{n}$, $\phi^n$, $\mu^n$, $r^n$, $T^n$, $Q^n, R^n, K^{n}$, we compute $\boldsymbol{u}^{n+1}$, $p^{n+1}$, $\phi^{n+1}$, $\mu^{n+1}$, $r^{n+1}$, $T^{n+1}$, $Q^{n+1}$, $R^{n+1}, K^{n+1}$ by the following three steps.\par
\noindent Step 1. 
\begin{align}
	& \frac{\phi^{n+1}-\phi^n}{\delta t}+Q^{n} \nabla \cdot\left(\boldsymbol{u}^{n}\phi^{n}\right)=\nabla\cdot \left(M(\phi^{n})\nabla\mu^{n+1}\right), \label{CH1_fst}\\
	& \mu^{n+1}=\lambda\left(-\epsilon \Delta \phi^{n+1}+\frac{s}{\epsilon} \phi^{n+1}+\frac{r^{n}}{\epsilon}(f(\phi^{n})-s\phi^{n})\right), \label{CH2_fst}\\
	&r^{n+1}-r^n=\alpha\bigg(-\left(E_{0}^{n+1}-E_{0}^{n}\right)+\frac{r^{n}\lambda}{\epsilon}\int_{\Omega}(f(\phi^{n})-s\phi^{n})(\phi^{n+1}-\phi^n)d\boldsymbol{x}\bigg), \label{rn_fst}\\
	& \frac{Q^{n+1}-Q^n}{\delta t}=\alpha Q^{n}\int_{\Omega}\left(\nabla \cdot\left(\boldsymbol{u}^{n} \phi^{n}\right) \mu^{n+1}+\left(\phi^{n} \nabla \mu^{n}\right) \cdot \boldsymbol{u}^{n+1}\right) d \boldsymbol{x}\label{Qn_fst}.
\end{align}
Step 2. 
\begin{align}
	&K^{n}\left(\rho^{n+1} \frac{\boldsymbol{u}^{n+1}-\boldsymbol{u}^{n}}{\delta t}+\frac{1}{2} \frac{\rho^{n+1}-\rho^n}{\delta t} \boldsymbol{u}^{n+1}\right)-\nabla \cdot\left(\nu^{n+1} D\left(\boldsymbol{u}^{n+1}\right)\right)+R^{n}\nabla\left(p^n+\psi^{n}\right) \label{NS1_fst} \\
	&+Q^{n}\phi^{n}\nabla\mu^{n}+R^{n} \rho^{n+1}\left(\boldsymbol{u}^{n} \cdot \nabla\right) \boldsymbol{u}^{n}+\frac{1}{2} R^{n}\nabla \cdot\left(\rho^{n+1} \boldsymbol{u}^{n}\right) \boldsymbol{u}^{n}+R^{n}\boldsymbol{J}^{n} \cdot \nabla \boldsymbol{u}^{n}+\frac{1}{2} R^{n} \nabla \cdot \boldsymbol{J}^{n} \boldsymbol{u}^{n}=\rho^{n+1}\mathbf{g},  \nonumber 
 \end{align}
 \begin{align}
	\frac{R^{n+1}-R^n}{\delta t}&= \alpha R^{n}\int_{\Omega}\left(\rho^{n+1}\left(\boldsymbol{u}^{n}\cdot \nabla\right) \boldsymbol{u}^{n} \cdot \boldsymbol{u}^{n+1}+\frac{1}{2} \nabla \cdot\left(\rho^{n+1} \boldsymbol{u}^{n}\right)\boldsymbol{u}^{n} \cdot \boldsymbol{u}^{n+1}\right) d \boldsymbol{x} \label{Rn_fst}\\
	&+\alpha R^{n}\int_{\Omega}\bigg(\left(\boldsymbol{J}^{n} \cdot \nabla \boldsymbol{u}^{n}\right) \cdot \boldsymbol{u}^{n+1}+\frac{1}{2} \nabla \cdot \boldsymbol{J}^{n} \boldsymbol{u}^{n} \cdot \boldsymbol{u}^{n+1}+\nabla\left(p^n+\psi^{n}\right)\cdot\boldsymbol{u}^{n+1}\bigg) d \boldsymbol{x}, \nonumber\\
	\frac{K^{n+1}-K^n}{\delta t}&= \alpha\int_{\Omega}-\frac{\frac{1}{2}\rho^{n+1}\|\boldsymbol{u}^{n+1}\|^{2}-\frac{1}{2}\rho^{n}\|\boldsymbol{u}^{n}\|^{2}}{\delta t}+K^{n}\left(\rho^{n+1} \frac{\boldsymbol{u}^{n+1}-\boldsymbol{u}^{n}}{\delta t}+\frac{1}{2} \frac{\rho^{n+1}-\rho^n}{\delta t} \boldsymbol{u}^{n+1}\right)\cdot\boldsymbol{u}^{n+1} d\boldsymbol{x}. \label{Kn_fst}
\end{align}
where $\boldsymbol{J}^n=\boldsymbol{J}\left(\mu^n\right)$.\\
Step 3. 
 \begin{align}
	&\Delta\psi^{n+1}=T^{n}\frac{\chi}{2\delta t} \nabla \cdot \boldsymbol{u}^{n+1}, \label{NS2_fst}\\
  &p^{n+1}=\psi^{n+1}+p^{n}. \label{NS3_fst}
\end{align}
where $\chi=\min(\rho_{1},\rho_{2})$.
\begin{align}
	\frac{T^{n+1}-T^n}{\delta t}&=\alpha T^{n}\int_{\Omega}p^{n+1}\nabla \cdot \boldsymbol{u}^{n+1}d\boldsymbol{x}. \label{Tn_fst}
\end{align}
The boundary conditions are 
$$
\left.\boldsymbol{u}^{n+1}\right|_{\partial \Omega}=\mathbf{0},\left.\quad \partial_{\boldsymbol{n}} \phi^{n+1}\right|_{\partial \Omega}=\left.\partial_{\boldsymbol{n}} \mu^{n+1}\right|_{\partial \Omega}=\left.\partial_{\boldsymbol{n}} p^{n+1}\right|_{\partial \Omega}=0,
$$
which also implies $\partial_{\boldsymbol{n}} \psi^{n+1}|_{\partial \Omega}=0$.
Periodic boundary conditions can also be assumed for all the variables.\\
\noindent The other notations used in the scheme read as
\begin{align}
\left\{\begin{array}{l}
\hat{\phi}= \begin{cases}\phi, & |\phi|<1, \\
\operatorname{sign}(\phi), & |\phi|>1,\end{cases} \\
\rho^{n+1}=\frac{\rho_1-\rho_2}{2} \hat{\phi}^{n+1}+\frac{\rho_1+\rho_2}{2}, \quad \nu^{n+1}=\frac{\nu_1-\nu_2}{2} \hat{\phi}^{n+1}+\frac{\nu_1+\nu_2}{2} .
\end{array}\right. \label{phi_hat}
\end{align}
$\mathbf{Remarks.}$ 
\begin{itemize}
    \item[1.] In Step 1,  \eqref{CH1_fst} and \eqref{CH2_fst} form a linear system for solving $\phi^{n+1}$ and $\mu^{n+1}$ which is decoupled from other equations. Once $\phi^{n+1}$ and $\mu^{n+1}$ are solved,  $r^{n+1}$ and $Q^{n+1}$ are then updated from \eqref{rn_fst} and \eqref{Qn_fst}.
    \item[2.] In Step 2, \eqref{NS1_fst} is an elliptic equation for $\boldsymbol{u}^{n+1}$. $R^{n+1}, K^{n+1}$ are updated explicitly from \eqref{Rn_fst} and \eqref{Kn_fst}.
     \item[3.] In Step 3, $p^{n+1}$ is solved from from \eqref{NS2_fst} and \eqref{NS3_fst}. $T^{n+1}$ are updated explicitly from \eqref{Tn_fst}. 
    \item[4.] We set $\beta$ in \eqref{MCHNS_9_sed} to zero because the accuracy of pressure will not be affected in the first-order scheme. 
\end{itemize}

\subsubsection{The second-order scheme} \label{second-order}
We next show that the energy stable second order scheme can also be designed for  \eqref{MCHNS_1_sed}-\eqref{MCHNS_9_sed}.
Given $\boldsymbol{u}^{n-1}$, $p^{n-1}$, $\phi^{n-1}$, $\mu^{n-1}$, $r^{n-1}$, $T^{n-1}$, $Q^{n-1}, R^{n-1}, K^{n-1}$ ,$\boldsymbol{u}^n$, $p^{n}$, $\phi^n$, $\mu^n$, $r^n$, $T^n$, $Q^n, R^n, K^{n}$, we compute $\boldsymbol{u}^{n+1}$, $\phi^{n+1}$, $\mu^{n+1}$, $r^{n+1}$, $T^{n+1}$, $Q^{n+1}$, $R^{n+1}$ by the following three steps.  In below, the superscript $*$ means $(\cdot)^{*}=2(\cdot)^{n}-(\cdot)^{n-1}$.\\
\noindent Step 1. 
\begin{align}
	& \frac{3\phi^{n+1}-4\phi^n+\phi^{n-1}}{2\delta t}+Q^{*} \nabla \cdot\left(\boldsymbol{u}^{*} \phi^{*}\right)=\nabla\cdot \left(M(\phi^{*})\nabla\mu^{n+1}\right), \label{CH1_sed}\\
	& \mu^{n+1}=\lambda\left(-\epsilon \Delta \phi^{n+1}+\frac{s}{\epsilon} \phi^{n+1}+\frac{r^{*}}{\epsilon}(f(\phi^{*})-s\phi^{*})\right), \label{CH2_sed}\\
	&3r^{n+1}-4r^n+r^{n-1}=\alpha\bigg(-\left(3E_{0}^{n+1}-4E_{0}^{n}+E_{0}^{n-1}\right) \nonumber \\
	&+\frac{r^{*}\lambda}{\epsilon}\int_{\Omega}(f(\phi^{*})-s\phi^{*})(3\phi^{n+1}-4\phi^n+\phi^{n-1})d\boldsymbol{x}\bigg), \label{rn_sed}\\
	& \frac{3Q^{n+1}-4Q^n+Q^{n-1}}{2\delta t}=\alpha Q^{*}\int_{\Omega}\left(\nabla \cdot\left(\boldsymbol{u}^{*} \phi^{*}\right) \mu^{n+1}+\left(\phi^{*} \nabla \mu^{*}\right) \cdot \boldsymbol{u}^{n+1}\right) d \boldsymbol{x}\label{Qn_sed}.
\end{align}
Step 2. 
\begin{align}
	&K^{*}\left(\rho^{n+1} \frac{3\boldsymbol{u}^{n+1}-4\boldsymbol{u}^{n}+\boldsymbol{u}^{n-1}}{2\delta t}+\frac{1}{2} \frac{3\rho^{n+1}-4\rho^n+\rho^{n-1}}{2\delta t} \boldsymbol{u}^{n+1}\right)-\nabla \cdot\left(\nu^{n+1} D\left(\boldsymbol{u}^{n+1}\right)\right)+R^{*}\nabla\left(p^n+\frac{4}{3}\omega^{n}-\frac{1}{3}\omega^{n-1}\right) \label{NS1_sed} \nonumber\\
	&+Q^{*}\phi^{*}\nabla\mu^{*}+R^{*} \rho^{n+1}\left(\boldsymbol{u}^{*} \cdot \nabla\right) \boldsymbol{u}^{*}+\frac{1}{2} R^{*} \nabla \cdot\left(\rho^{n+1} \boldsymbol{u}^{*}\right) \boldsymbol{u}^{*}+R^{*} \boldsymbol{J}^{*} \cdot \nabla \boldsymbol{u}^{*}+\frac{1}{2} R^{*} \nabla \cdot \boldsymbol{J}^{*} \boldsymbol{u}^{*}=\rho^{n+1}\mathbf{g}, 
\end{align}
\begin{align}
	\frac{3R^{n+1}-4R^n+R^{n-1}}{2\delta t}&= \alpha R^{*}\int_{\Omega}\left(\rho^{n+1}\left(\boldsymbol{u}^{*} \cdot \nabla\right) \boldsymbol{u}^{*} \cdot \boldsymbol{u}^{n+1}+\frac{1}{2} \nabla \cdot\left(\rho^{n+1} \boldsymbol{u}^{*}\right) \boldsymbol{u}^{*} \cdot \boldsymbol{u}^{n+1}\right) d \boldsymbol{x} \label{Rn_sed}\\
	&+\alpha R^{*}\int_{\Omega}\left(\left(\boldsymbol{J}^{*} \cdot \nabla \boldsymbol{u}^{*}\right) \cdot \boldsymbol{u}^{n+1}+\frac{1}{2} \nabla \cdot \boldsymbol{J}^{*} \boldsymbol{u}^{*} \cdot \boldsymbol{u}^{n+1}+\nabla\left(p^n+\frac{4}{3}\omega^{n}-\frac{1}{3}\omega^{n-1}\right)\cdot\boldsymbol{u}^{n+1}\right) d \boldsymbol{x}, \nonumber \\
	\frac{3K^{n+1}-4K^n+K^{n-1}}{2\delta t}&= \alpha\int_{\Omega}-\frac{3\frac{1}{2}\rho^{n+1}\|\boldsymbol{u}^{n+1}\|^{2}-4\frac{1}{2}\rho^{n}\|\boldsymbol{u}^{n}\|^{2}+\frac{1}{2}\rho^{n-1}\|\boldsymbol{u}^{n-1}\|^{2}}{2\delta t}\nonumber \\
	&+K^{*}\left(\rho^{n+1} \frac{3\boldsymbol{u}^{n+1}-4\boldsymbol{u}^{n}+\boldsymbol{u}^{n-1}}{2\delta t}+\frac{1}{2} \frac{3\rho^{n+1}-4\rho^n+\rho^{n-1}}{2\delta t} \boldsymbol{u}^{n+1}\right)\cdot\boldsymbol{u}^{n+1} d\boldsymbol{x}. \label{Kn_sed}
\end{align} 
where $\rho^{n+1}$ and $\eta^{n+1}$ are calculated using \eqref{phi_hat}.\\
Step 3. 
\begin{align}
	&\Delta\omega ^{n+1}=T^{*}\frac{3\chi}{2\delta t} \nabla \cdot \boldsymbol{u}^{n+1},\label{NS2_sed}\\ 
 &p^{n+1}  =\omega^{n+1}+p^{n}-T^{*}\nu^{n+1}\nabla\cdot\boldsymbol{u}^{n+1},
\end{align}
\begin{align}
	\frac{3T^{n+1}-4T^n+T^{n-1}}{2\delta t}&=\alpha T^{*}\left(\int_{\Omega}p^{n+1}\nabla \cdot \boldsymbol{u}^{n+1}+\beta\nabla\left(\nu^{n+1}\nabla\cdot\boldsymbol{u}^{n+1}\right)\cdot\nabla p^{n+1} d\boldsymbol{x}\right), \label{Tn_sed}
\end{align}
 where $\chi=\min(\rho_{1},\rho_{2})$ and $\beta=\frac{2\delta t}{3\chi}$.
The boundary conditions read as
$$
\left.\boldsymbol{u}^{n+1}\right|_{\partial \Omega}=\mathbf{0},\left.\quad \partial_{\boldsymbol{n}} \phi^{n+1}\right|_{\partial \Omega}=\left.\partial_{\boldsymbol{n}} \mu^{n+1}\right|_{\partial \Omega}=\left.\partial_{\boldsymbol{n}} \omega^{n+1}\right|_{\partial \Omega}=0.
$$
Periodic boundary conditions can also be assumed for all the variables.\par
$\mathbf{Remarks.}$
\begin{itemize}
    \item[1.]
    In equations \eqref{CH1_sed}-\eqref{Tn_sed}, only linear systems need to be solved for ($\phi$, $\mu$) and ($\boldsymbol{u}$, $p$). The auxiliary variables ($Q$, $R$, $K$, $T$) are explicitly updated at each time step. Consequently, our method demonstrates higher computational efficiency compared to the DSAV method, which is verified through numerical simulations in Section \ref{ex_1}.

\item[2.] Within equation \eqref{NS2_sed}, the term $\nu^{n+1}\nabla\cdot\boldsymbol{u}^{n+1}$ is introduced to enhance the accuracy of pressure. The standard penalty method is afflicted by non-physical pressure boundary conditions, leading to a numerical boundary layer that ultimately reduces the precision of the pressure approximation \cite{guermond2004error}. In \cite{timmermans1996approximate}, a modified scheme, known as the rotation penalty method, is proposed to tackle this issue by incorporating the term $\nu^{n+1}\nabla\cdot\boldsymbol{u}^{n+1}$ in the pressure Poisson equation. For more detailed explanation, consult references \cite{guermond2006overview, guermond2004error}. The introduction of the term $\nu^{n+1}\nabla\cdot\boldsymbol{u}^{n+1}$ also presents challenges in establishing energy stability. Leveraging its property of zero energy contribution, we address this challenge by incorporating $\beta\nabla\left(\nu^{n+1}\nabla\cdot\boldsymbol{u}^{n+1}\right)\cdot\nabla p^{n+1}$ into equation \eqref{Tn_sed}.

\item[3.]  The modified rotational penalty method is employed in the aforementioned scheme to discretize the Navier-Stokes equations. The computation of pressure is decoupled from the momentum equation. This method can also avoid solving the elliptic equation with $\frac{1}{\rho}$ as the variable coefficients so that only one pressure Poisson equation with constant coefficients needs to be solved, which can greatly reduce the computational cost. In contrast to the standard penalty method, the modified rotational penalty method attains superior accuracy of pressure, and the computational overhead rises mildly by solving an extra scalar equation \eqref{Tn_sed} whose computational burden is negligible. {The artificial compression method can also be used to establish a numerical scheme that decouples pressure and velocity. However, it has the drawback of introducing non-physical acoustic waves (refer to \cite{wang2024efficient}). The modified rotational penalty method does not have this issue (see Section \ref{bubble_rising2D}).}

\item[4.]  It is evident that in numerical computations, $r^{n+1}$, $Q^{n+1}$, $R^{n+1}$, $T^{n+1}$, and $K^{n+1}$ generally do not preserve an exact value of 1. Hence, we introduce a parameter $\alpha$ on the right side of equations \eqref{rn_sed}, \eqref{Qn_sed}, \eqref{Kn_sed}, \eqref{Rn_sed}, and \eqref{Tn_sed}. Typically, $\alpha$ is chosen to be very small to ensure that the right-hand side remains sufficiently diminutive. Our numerical findings indicate that by selecting a suitably small $\alpha$, the values of $r^{n+1}$, $Q^{n+1}$, $R^{n+1}$, $T^{n+1}$, and $K^{n+1}$ are kept close to 1.

{ \item[5.]  In equations \eqref{rn_sed} and \eqref{Kn_sed}, a relationship is established between the conservative numerical time derivative and the non-conservative numerical time derivative of the free energy $E_{0}$ and the kinetic energy $\frac{1}{2}\rho\boldsymbol{u}^{2}$. In the energy proof, by merging \eqref{CH2_sed} with \eqref{rn_sed} (or \eqref{NS1_sed} with \eqref{Kn_sed}), the non-conservative numerical time derivative (the third term on the right-hand side of \eqref{CH2_sed}) will be substituted with its conservative form.}

\end{itemize}

Energy stable property can be proved for both the first order scheme \eqref{CH1_fst}-\eqref{Tn_fst} and the second order scheme \eqref{CH1_sed}-\eqref{Tn_sed}. We here prove the following Theorem for the second-order scheme. \\
Theorem 5. The solutions of the time-discrete scheme \eqref{CH1_sed}-\eqref{Tn_sed} satisfy the following energy law
\begin{align}
	E^{n+1} \leqslant E^n-\frac{\delta t}{2}\left\|\sqrt{\nu^{n+1}} D\left(\boldsymbol{u}^{n+1}\right)\right\|^2-\delta t\int_{\Omega}M(\phi^{*})|\nabla \mu^{n+1}|^2 d \boldsymbol{x}, \label{EM_de_sed}
\end{align}
where $\sigma^k=\sqrt{\rho^k}$ for any $k$ and
\begin{align}
	E^{n+1}= &\frac{1}{2}\left(\frac{3}{2}\|\sigma^{n+1}\boldsymbol{u}^{n+1}\|^{2}-\frac{1}{2}\|\sigma^{n}\boldsymbol{u}^{n}\|^{2}\right)+\frac{\lambda\epsilon}{4}\left(\|\nabla\phi^{n+1}\|^{2}+\|2\nabla\phi^{n+1}-\nabla\phi^{n}\|^{2}\right) \nonumber \\
	&+\frac{\lambda s}{4\epsilon}\left(\|\phi^{n+1}\|^{2}+\|2\phi^{n+1}-\phi^{n}\|^{2}\right)+\frac{\delta t^{2}}{3\chi}\|\nabla p^{n+1}\|^{2}+\frac{1}{2\alpha}\left(3r^{n+1}-r^{n}\right)+\frac{1}{2}\left(3E_{0}^{n+1}-E_{0}^{n}\right) \nonumber\\
	&+\frac{1}{2\alpha}\left(3R^{n+1}-R^{n}\right)+\frac{1}{2\alpha}\left(3Q^{n+1}-Q^{n}\right)+\frac{1}{2\alpha}\left(3T^{n+1}-T^{n}\right)+\frac{1}{2\alpha}\left(3K^{n+1}-K^{n}\right). \label{EM_sed}
\end{align}
\textbf{Proof.} By taking the inner product of \eqref{NS1_sed} with $2\delta t\boldsymbol{u}^{n+1}$ in the $L^2$ space, we obtain
\begin{align}
	&K^{*}\left(\rho^{n+1},\left(3\boldsymbol{u}^{n+1}-4\boldsymbol{u}^{n}+\boldsymbol{u}^{n-1}\right)\cdot\boldsymbol{u}^{n+1}\right)+\frac{1}{2}\left(3\rho^{n+1}-4\rho^{n}+\rho^{n-1}, \boldsymbol{u}^{n+1}\cdot\boldsymbol{u}^{n+1}\right) \nonumber\\
	&+\delta t\left\|\sqrt{\nu^{n+1}} D\left(\boldsymbol{u}^{n+1}\right)\right\|^2+2\delta tR^{*}\int_{\Omega}\nabla\left(p^{n}+\frac{4}{3}\omega^{n}-\frac{1}{3}\omega^{n-1}\right)\cdot \boldsymbol{u}^{n+1}d\boldsymbol{x}\nonumber \\
	&+2\delta tQ^{*}\int_{\Omega}\phi^{*}\nabla\mu^{*}\cdot\boldsymbol{u}^{n+1} d\boldsymbol{x}+2\delta t R^{*}\int_{\Omega}\rho^{n+1}\left(\boldsymbol{u}^{*}\cdot\nabla\right)\boldsymbol{u}^{*}\cdot\boldsymbol{u}^{n+1}d\boldsymbol{x}\nonumber \\
	&+\delta tR^{*}\int_{\Omega}\nabla \cdot (\rho^{n+1}\boldsymbol{u}^{*})\boldsymbol{u}^{*}\cdot\boldsymbol{u}^{n+1}d\boldsymbol{x}+2\delta tR^{*}\int_{\Omega}\left(J^{*}\cdot\nabla\boldsymbol{u}^{*}\right)\cdot \boldsymbol{u}^{n+1}d\boldsymbol{x}\nonumber\\
	&+\delta t R^{*}\int_{\Omega}\nabla\cdot J^{*}\boldsymbol{u}^{*}\cdot\boldsymbol{u}^{n+1} d\boldsymbol{x}=0.\label{NS1_v1_sed}
\end{align}
By combining \eqref{Rn_sed}, \eqref{Kn_sed} and \eqref{NS1_v1_sed}, we derive
\begin{align}
	&\int_{\Omega}3\frac{1}{2}\rho^{n+1}|\boldsymbol{u}^{n+1}|^{2}-4\frac{1}{2}\rho^{n}|\boldsymbol{u}^{n}|^{2}+\frac{1}{2}\rho^{n-1}|\boldsymbol{u}^{n-1}|^{2}d\boldsymbol{x}+\frac{1}{\alpha}\left(3K^{n+1}-4K^{n}+K^{n-1}\right)\nonumber\\
	&+\frac{1}{\alpha}\left(3R^{n+1}-4R^{n}+R^{n-1}\right)=-\delta t\left\|\sqrt{\nu^{n+1}} D\left(\boldsymbol{u}^{n+1}\right)\right\|^2-2\delta tQ^{*}\int_{\Omega}\phi^{*}\nabla\mu^{*}\cdot\boldsymbol{u}^{n+1} d\boldsymbol{x}. \label{NS2_v3_sed}
\end{align}
We take the $L^2$ inner product of \eqref{CH1_sed} with $4\delta t \mu^{n+1}$ and use integration by parts to get
\begin{align}
	2\left(\mu^{n+1},3\phi^{n+1}-4\phi^{n}+\phi^{n-1}\right)+4\delta t Q^{*}\left(\nabla\cdot (\boldsymbol{u}^{*}\phi^{*}),\mu^{n+1}\right)=-4\delta t \int_{\Omega}M(\phi^{*})|\nabla \mu^{n+1}|^2 d \boldsymbol{x}. \label{CH1_v1_sed}
\end{align}
By taking the $L^2$ inner product of \eqref{CH2_sed} with $2\left(3\phi^{n+1}-4\phi^n+\phi^{n-1}\right)$ and using integration by parts and using the following identity:
\begin{align}
2a(3a-4b+c)=a^2-b^2+(2a-b)^2-(2b-c)^2+(a-2b+c)^2, \label{Second_order}
\end{align}we get
\begin{align}
	&2\left(\mu^{n+1},3\phi^{n+1}-4\phi^{n}+\phi^{n-1}\right)=\lambda\epsilon\big(\|\nabla\phi^{n+1}\|^{2}-\|\nabla\phi^{n}\|^{2}+\|2\nabla\phi^{n+1}-\nabla\phi^{n}\|^{2}-\|2\nabla\phi^{n}-\nabla\phi^{n-1}\|^{2}\nonumber \\
	&+\|\nabla\phi^{n+1}-2\nabla\phi^{n}+\nabla\phi^{n-1}\|^{2}\big)+\frac{2\lambda r^{*}}{\epsilon}\left(f(\phi^{*})-s\phi^{*},3\phi^{n+1}-4\phi^{n}+\phi^{n-1}\right)\nonumber \\
	&+\frac{\lambda s}{\epsilon}\left(\|\phi^{n+1}\|^{2}-\|\phi^{n}\|^{2}+\|2\phi^{n+1}-\phi^{n}\|^{2}-\|2\phi^{n}-\phi^{n-1}\|^{2}+\|\phi^{n+1}-2\phi^{n}+\phi^{n-1}\|^{2}\right).\label{CH2_v1_sed}
\end{align}
By combining \eqref{CH1_v1_sed} and \eqref{CH2_v1_sed},
\begin{align}
&\frac{\lambda\epsilon}{2}\left(\|\nabla\phi^{n+1}\|^{2}-\|\nabla\phi^{n}\|^{2}+\|2\nabla\phi^{n+1}-\nabla\phi^{n}\|^{2}-\|2\nabla\phi^{n}-\nabla\phi^{n-1}\|^{2}+\|\nabla\phi^{n+1}-2\nabla\phi^{n}+\nabla\phi^{n-1}\|^{2}\right)\nonumber \\
&+\frac{\lambda s}{2\epsilon}\left(\|\phi^{n+1}\|^{2}-\|\phi^{n}\|^{2}+\|2\phi^{n+1}-\phi^{n}\|^{2}-\|2\phi^{n}-\phi^{n-1}\|^{2}+\|\phi^{n+1}-2\phi^{n}+\phi^{n-1}\|^{2}\right) \nonumber \\
&+\frac{\lambda r^{*}}{\epsilon}\left(f(\phi^{*})-s\phi^{*},3\phi^{n+1}-4\phi^{n}+\phi^{n-1}\right)=-2\delta t \int_{\Omega}M(\phi^{*})|\nabla \mu^{n+1}|^2 d \boldsymbol{x}-2\delta tQ^{*}\left(\nabla\cdot (\boldsymbol{u}^{*}\phi^{*}),\mu^{n+1}\right). \label{CH_v1_sed}
\end{align}
By combining \eqref{CH_v1_sed} and \eqref{rn_sed},
\begin{align}
	&\frac{\lambda\epsilon}{2}\left(\|\nabla\phi^{n+1}\|^{2}-\|\nabla\phi^{n}\|^{2}+\|2\nabla\phi^{n+1}-\nabla\phi^{n}\|^{2}-\|2\nabla\phi^{n}-\nabla\phi^{n-1}\|^{2} + \|\nabla\phi^{n+1}-2\nabla\phi^{n}+\nabla\phi^{n-1}\|^{2}\right) \nonumber \\
	&+\frac{\lambda s}{2\epsilon}\left(\|\phi^{n+1}\|^{2}-\|\phi^{n}\|^{2}+\|2\phi^{n+1}-\phi^{n}\|^{2}-\|2\phi^{n}-\phi^{n-1}\|^{2}+\|\phi^{n+1}-2\phi^{n}+\phi^{n-1}\|^{2}\right) \nonumber \\
	&+\frac{1}{\alpha}(3r^{n+1}-4r^n+r^{n-1})+\left(3E_{0}^{n+1}-4E_{0}^{n}+E_{0}^{n-1}\right)=-2\delta t \int_{\Omega}M(\phi^{*})|\nabla \mu^{n+1}|^2 d \boldsymbol{x}-2\delta tQ^{*}\left(\nabla\cdot (\boldsymbol{u}^{*}\phi^{*}),\mu^{n+1}\right).  \label{CH_v2_sed}
\end{align}
We take the $L^2$ inner product of \eqref{NS2_sed} with $4\frac{\delta t^{2}}{3\chi}p^{n+1}$ and use integration by parts to get
\begin{align}
	\frac{2\delta t^{2}}{3\chi}\left(\|\nabla p^{n+1}\|^{2}-\|\nabla p^{n}\|^{2}+\|\nabla p^{n+1}-\nabla p^{n}\|^{2}\right)&=-2\delta tT^{*}\left(\nabla\cdot\boldsymbol{u}^{n+1},p^{n+1}\right) \nonumber\\
 &-\frac{4\delta t^{2}}{3\chi}T^{*}\int_{\Omega}\nabla\left(\nu^{n+1}\nabla\cdot\boldsymbol{u}^{n+1}\right)\cdot\nabla p^{n+1} d\boldsymbol{x}.\label{p_v1_sed}
\end{align}
By combining  \eqref{Qn_sed}, \eqref{Rn_sed}, \eqref{Tn_sed}, \eqref{NS2_v3_sed}, \eqref{CH_v2_sed} and  \eqref{p_v1_sed}, we derive
\begin{align}
	&\frac{1}{2}\left(3\|\sigma^{n+1}\boldsymbol{u}^{n+1}\|^{2}-4\|\sigma^{n}\boldsymbol{u}^{n}\|^{2}+\|\sigma^{n-1}\boldsymbol{u}^{n-1}\|^{2}\right)+\frac{1}{\alpha}\left(3K^{n+1}-4K^{n}+K^{n-1}\right)+\frac{1}{\alpha}\left(3Q^{n+1}-4Q^{n}+Q^{n-1}\right)\nonumber\\
	&\frac{\lambda\epsilon}{2}\left(\|\nabla\phi^{n+1}\|^{2}-\|\nabla\phi^{n}\|^{2}+\|2\nabla\phi^{n+1}-\nabla\phi^{n}\|^{2}-\|2\nabla\phi^{n}-\nabla\phi^{n-1}\|^{2} + \|\nabla\phi^{n+1}-2\nabla\phi^{n}+\nabla\phi^{n-1}\|^{2}\right) \nonumber \\
	&+\frac{\lambda s}{2\epsilon}\left(\|\phi^{n+1}\|^{2}-\|\phi^{n}\|^{2}+\|2\phi^{n+1}-\phi^{n}\|^{2}-\|2\phi^{n}-\phi^{n-1}\|^{2}+\|\phi^{n+1}-2\phi^{n}+\phi^{n-1}\|^{2}\right) \nonumber \\
	&+\frac{1}{\alpha}(3r^{n+1}-4r^n+r^{n-1})+\frac{1}{\alpha}(3T^{n+1}-4T^n+T^{n-1})+\frac{1}{\alpha}\left(3R^{n+1}-4R^{n}+R^{n-1}\right)+\left(3E_{0}^{n+1}-4E_{0}^{n}+E_{0}^{n-1}\right)\nonumber \\
	&+\frac{2\delta t^{2}}{3\chi}\left(\|\nabla p^{n+1}\|^{2}-\|\nabla p^{n}\|^{2}+\|\nabla p^{n+1}-\nabla p^{n}\|^{2}\right)=-\delta t\left\|\sqrt{\nu^{n+1}} D\left(\boldsymbol{u}^{n+1}\right)\right\|^2-2\delta t \int_{\Omega}M(\phi^{*})|\nabla \mu^{n+1}|^2 d \boldsymbol{x}. \label{dEM_decay_sed}
\end{align}
Finally, we obtain \eqref{EM_de_sed} from \eqref{dEM_decay_sed} after dropping some unnecessary positive terms.\\
$\mathbf{Remark.} $ For the proof of the first-order-in-time numerical scheme, it is only necessary to replace \eqref{Second_order} with the following equation
\begin{equation}
2(a-b) a=a^2-b^2+(a-b)^2. \label{first-order}
\end{equation}
{$\mathbf{Remark.} $ When accounting for the impact of gravity force, it becomes necessary to incorporate an additional constant auxiliary variable to address the term $\rho\mathbf{g}$. However, since our current emphasis does not lie on this particular term, we will not discuss it in detail.}
\subsection{Numerical results}
In this section, we numerically investigate the accuracy, energy stability, and efficiency of the proposed schemes \eqref{CH1_fst}-\eqref{Tn_fst} and \eqref{CH1_sed}-\eqref{Tn_sed}. We present numerical examples for various benchmark scenarios, including the ascent of lighter droplets under gravity in both 2D and 3D settings, and the Rayleigh-Taylor instability in 2D. The computational domain is configured as rectangular domains in either 2D or 3D. Spatial discretization is accomplished by utilizing the finite element method. 

\subsubsection{Accuracy and stability test} \label{ex_1}
\begin{figure}[H]
\centering
	\includegraphics[width=75mm]{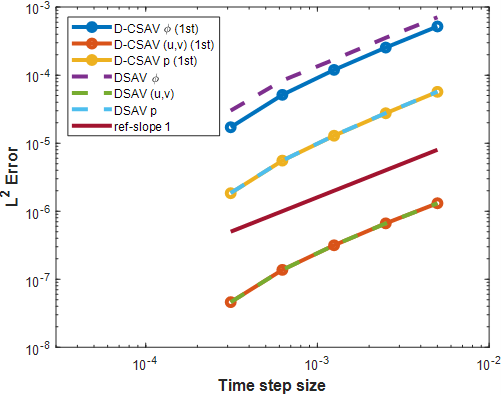}
\caption{First-order temporal convergence tests: $L^{2}$ errors at $t=0.64$ for the phase-field variable $\phi$, the average of the two velocity components $(u,v)$, and the pressure $p$.}
\label{fig:Fig.1all}
\end{figure}

\begin{figure}[H]
\centering
\subfigure[]{
	\label{fig:Fig 2}
	\includegraphics[width=75mm]{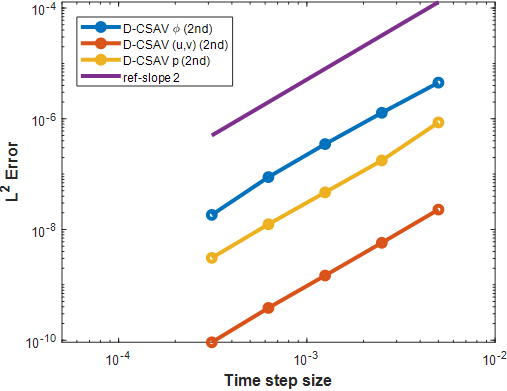}}
\subfigure[]{
	\label{fig:Fig 3}
	\includegraphics[width=75mm]{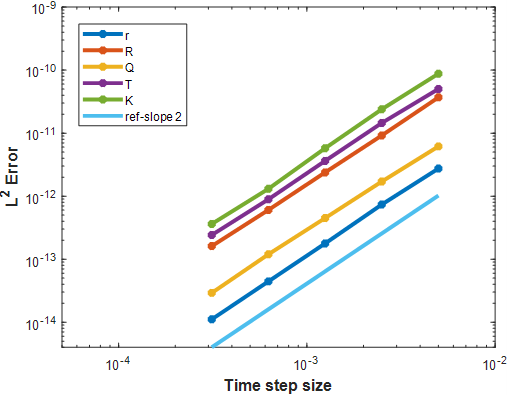}}
\caption{Second-order temporal convergence tests: 
(a) $L^{2}$ errors of the phase-field variable $\phi$, the average of the two velocity components $(u,v)$, and the pressure $p$ at $t=0.64$.
(b) $L^{2}$ errors of the scalar variables $r, Q, R, T, K$ at $t=0.64$.}
\label{fig:Fig.2all}
\end{figure}

\begin{figure}[H]
\centering
\subfigure[]{
	\label{fig:Fig3_1}
\includegraphics[width=75mm]{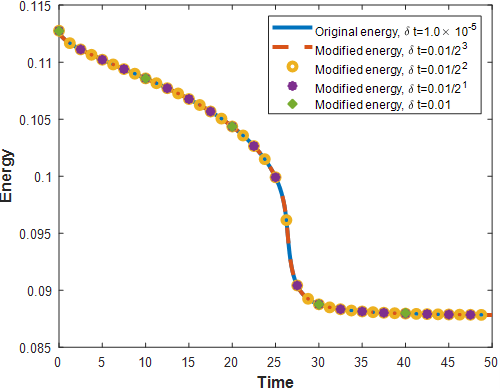}}
\subfigure[]{
	\label{fig:Fig3_2}
\includegraphics[width=75mm]{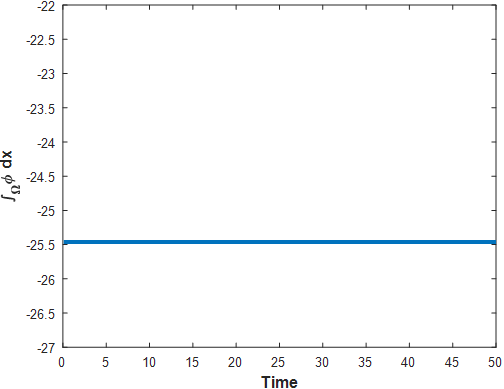}
 }
\caption{(a) Evolution of the original energy \eqref{EO} computed using $\delta t=10^{-5}$ and the modified total free energy \eqref{EM_sed} computed using different time step sizes. (b)
Evolution of the volume $\int_{\Omega}\phi dx$.
}
\label{fig:Fig.3all}
\end{figure}

We set the 2D computational domain to be $\Omega=[0,2 \pi]^2$. The system is assumed to satisfy the periodic boundary conditions. The initial conditions are set as follows:
$$
\left\{\begin{array}{l}
	\phi^0(x, y)=1+\sum_{i=1}^2 \tanh \left(\frac{r_i-\sqrt{\left(x-x_i\right)^2+\left(y-y_i\right)^2}}{1.5 \epsilon}\right), \\
	u^0(x, y)=0, \quad p^0(x, y)=0,
\end{array}\right.
$$
where
$$
r_1=1.4, \quad r_2=0.5, \quad x_1=\pi-0.8, \quad x_2=\pi+1.7, \quad y_1=y_2=\pi.
$$
The model parameters are given as
$$
\lambda=0.01, \quad \epsilon=0.08, \quad M(\phi)=1,\quad \alpha=10^{-2}, \quad\rho_1=10,\quad \rho_2=1, \quad\nu_1=1, \quad\nu_2=1, \quad s=4.
$$
To verify the temporal convergence, a sufficiently fine mesh of size $256\times 256$ is used to resolve the interface. This ensures that the spatial discretization error is negligible compared to the temporal discretization error. Due to the absence of an exact solution for this problem, we take the numerical solution obtained using a small mesh of size $512\times512$ and a small time step size $\delta t=10^{-6}$ as the reference. In Figure \ref{fig:Fig.1all}, we present the $L^2$ errors of $\phi$, the average of the two
velocity components $\frac{u+v}{2}$, and $p$ calculated at $t=0.64$ using a first-order-in-time D-CSAV scheme \eqref{CH1_fst}-\eqref{Tn_fst} and the DSAV method \cite{chen2022highly}. It can be seen that a good first-order temporal convergence order is achieved for all variables in both methods, while D-CSAV method results in a slightly better accuracy for $\phi$. In terms of the computational efficiency, we show in Table \ref{table1} a comparison of compute time using DSAV and D-CSAV for solving the problem up to $t=0.64$, both with first-order temporal accuracy. Notably, D-CSAV is more efficient as it requires only half of the compute time for DSAV method. {For the second-order-in-time D-CSAV scheme, we show in Figure \ref{fig:Fig.2all} the $L^2$ errors of  $\phi$, $\frac{u+v}{2}$, $p$ and $r, Q, R, T, K$ at $t=0.64$ (note that the exact solution of scalar auxiliary variables is 1).} The results show a good second-order temporal convergence rate. With respect to the pressure $p$, a superconvergence effect is probable owing to the regularity of the domain. A similar superconvergence effect was reported for the rotational version of the pressure-correction algorithm for constant density flows \cite{guermond2004error}.
In  Figure \ref{fig:Fig3_1}, we plot the evolution of the original energy \eqref{EO} computed using $\delta t=10^{-5}$ and the modified total free energy \eqref{EM_sed} computed using different time steps (after subtracting some constant auxiliary variables as mentioned in \eqref{EOEM}). 
All modified energy curves are monotonically decreasing and consistent with the original energy curve, illustrating the unconditional stability of the method. Figure \ref{fig:Fig3_2} illustrates the time evolution of the volume $\int_{\Omega}\phi dx$, which remains almost constant, demonstrating excellent volume conservation property.


\begin{table}[h!]
\centering
\begin{tabular}{rrr}
\hline & DSAV & D-CSAV \\
\hline$\delta t$  & CPU time(s)  & CPU time(s) \\
\hline $3.125 \times 10^{-4}$ & $9.867 \times 10^3$ & $5.472 \times 10^3$ \\
$6.25  \times 10^{-4}$ &  $5.050\times 10^3$ & $2.684 \times 10^3$ \\
$1.25 \times 10^{-3}$ &  $2.574 \times 10^3$ & $1.258 \times 10^3$ \\
$2.5  \times 10^{-3}$ &  $1.298 \times 10^3$  & $6.382 \times 10^2$ \\
$5  \times 10^{-3}$ & $6.417 \times 10^2$  & $3.214 \times 10^2$ \\
\hline
\end{tabular}
\caption{A comparison of compute time between DSAV and D-CSAV, both with first-order temporal accuracy.}
\label{table1}
\end{table}


\subsubsection{Bubble rising} \label{bubble_rising2D}
To further validate the accuracy of the method, we consider two benchmark cases for a single bubble rising in a calm water channel. In both  cases, the domain $[0,1] \times[0,2]$ consists of fluid 1 $(\phi \approx 1)$ except for a circular bubble filled with fluid 2 $(\phi \approx-1)$.
The gravitational source term $\mathbf{g}=(0,-0.98)$ is used in the momentum equation. The initial condition for $\phi$ is taken as
$$
\phi_0(x, y)=\tanh \left(\frac{r-0.25}{\sqrt{2} \epsilon}\right), \quad r=\sqrt{(x-0.5)^2+(y-0.5)^2}.
$$
We use a homogeneous Neumann boundary condition for the phase-field variables $\phi$ and $\mu$, and no-slip velocity boundary condition on the horizontal boundaries and free slip velocity boundary condition on the vertical boundaries. Referring to studies by Hosseini et al. \cite{hosseini2017isogeometric} and Fu et al. \cite{fu2021linear}, in Case 1, the bubble attains a stable ellipsoidal shape; in Case 2, it acquires a non-convex form with thin filaments that are prone to breaking, attributed to reduced surface tension and a higher density ratio \cite{hosseini2017isogeometric}.\par
Following \cite{hosseini2017isogeometric}, we consider the following physical parameters for the two cases, respectively,\\
Case 1: $\rho_1=1000$, $\rho_2=100$, $\eta_1=10$, $\eta_2=1$, $\sigma=24.5$, $\lambda=\frac{3}{2\sqrt{2}}\sigma$, $M(\phi)=\gamma(\phi^{2}-1)^{2}$, $\gamma=0.00004$, $\epsilon=0.01$, $\Delta x=2^{-6}$, $\delta t= 10^{-4}$, $\alpha=10^{-5}$.\\
Case 2: $\rho_1=1000$, $\rho_2=1$, $\eta_1=10$, $\eta_2=0.1$, $\sigma=1.96$, $\lambda=\frac{3}{2\sqrt{2}}\sigma$, $M(\phi)=\gamma(\phi^{2}-1)^{2}$, $\gamma=0.00026$, $\epsilon=0.01$, $\Delta x=2^{-6}$, $\delta t=0.5\times 10^{-4}$, $\alpha=10^{-5}$.\par
\begin{figure}[H]
\centering
\subfigure[]{
	\label{bubble motion}
	\includegraphics[width=60mm]{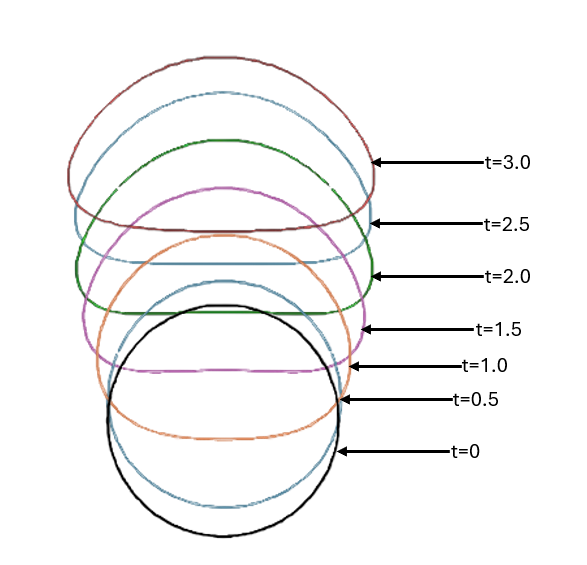}}
\subfigure[]{
	\label{bubble shape}
	\includegraphics[width=40mm]{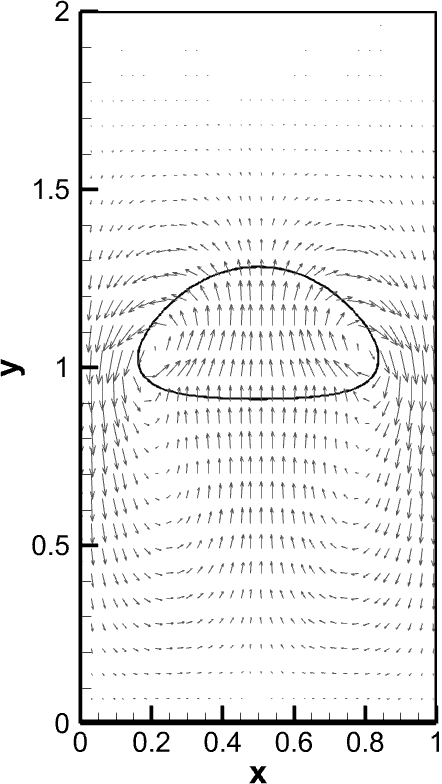}}
\caption{(a) Shapes of the rising bubble at different times for Case 1. (b) The shape and velocity vector distribution of the rising bubble at final time $t=3$ for Case 1.}
\label{b_Case1}
\end{figure}
\begin{figure}[H]
\centering
\includegraphics[width=150mm]{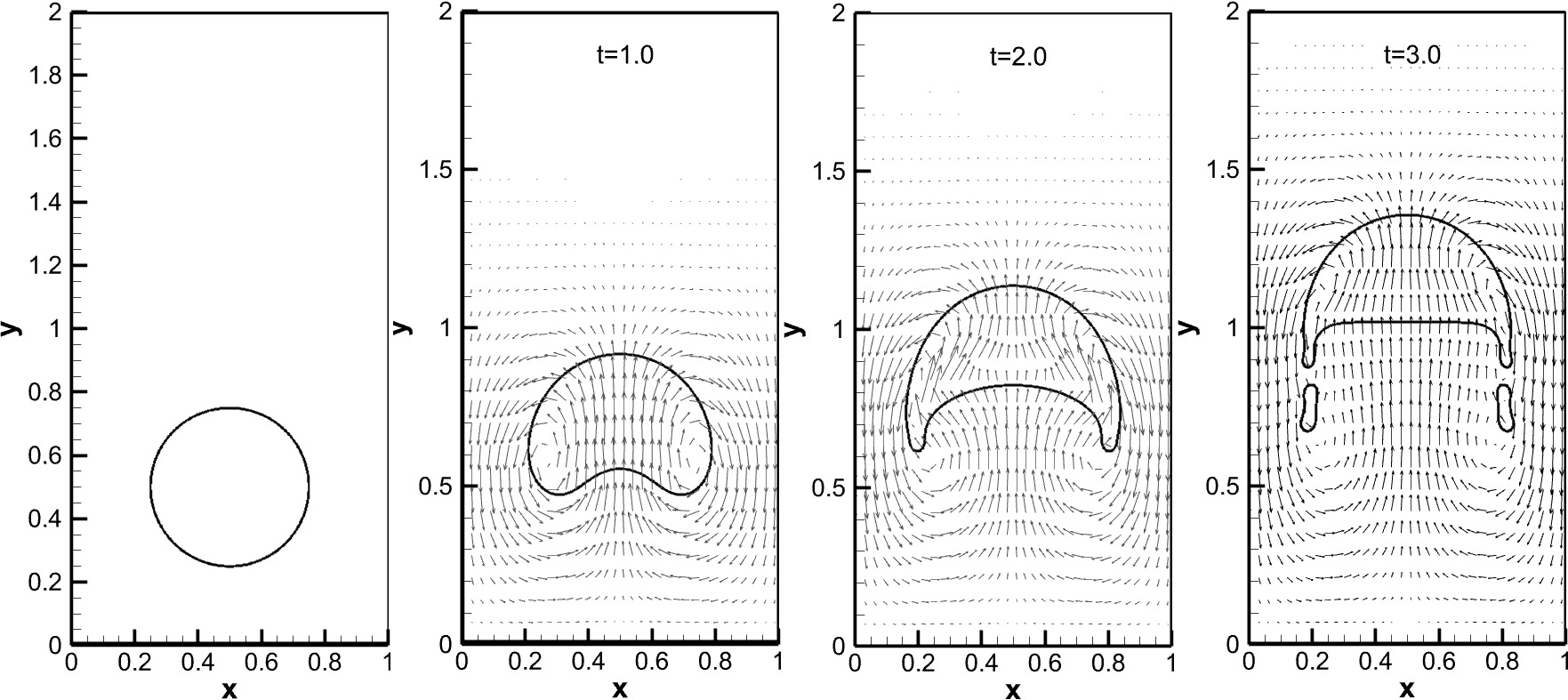}
\caption{Shapes of the rising bubble and the velocity vector distribution at different times for Case 2.}
\label{bubble motion_case2}
\end{figure}

Figure \ref{b_Case1} and \ref{bubble motion_case2} illustrate, at different time steps, the shapes of the rising bubble for Case 1 and Case 2 obtained using the proposed second-order-in-time numerical scheme \eqref{CH1_sed}-\eqref{Tn_sed}. Contour plots of the phase field variable $\phi$ and the velocity vector distribution diagram at the final time $t=3$ are provided in the figures. We observe the expected stable ellipsoidal bubble for Case 1 and the non-convex shape with filaments for Case 2. In particular, our simulation results in a break of the bubble filaments for Case 2.

Next, to validate the accuracy of the method, we compute the following benchmark quantities by using the first-order-in-time scheme \eqref{CH1_fst}-\eqref{Tn_fst} and the second-order-in-time scheme \eqref{CH1_sed}-\eqref{Tn_sed}.\\
Center of mass:
$$
y_c=\frac{\int_{\phi<0} y d\boldsymbol{x}}{\int_{\phi<0} 1 d\boldsymbol{x}},
$$
where $y$ is the vertical coordinate.\\
Rising velocity:
$$
V_c=\frac{\int_{\phi<0} v d\boldsymbol{x}}{\int_{\phi<0} 1 d\boldsymbol{x}},
$$
where $v$ is the vertical component of the velocity $\boldsymbol{u}$.
Then, we compare them with reference values from \cite{hysing2009quantitative}.


\begin{figure}[H]
\centering
\subfigure[]{
	\label{case1_Mass}
	\includegraphics[width=73mm]{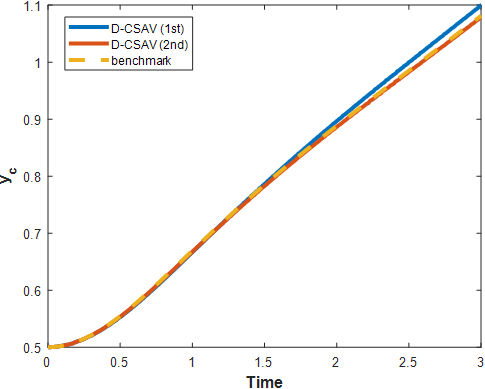}}
\subfigure[]{
	\label{case1_vel}
	\includegraphics[width=75mm]{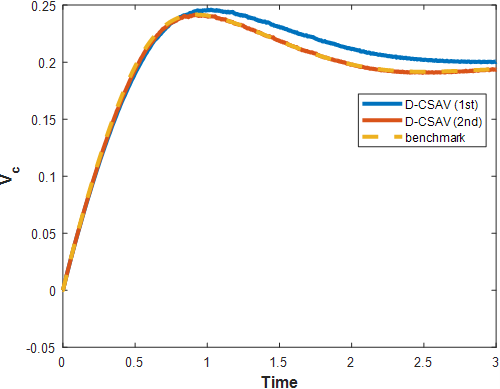}}
\caption{(a) Center of mass and (b) rising velocity obtained using the first and second order D-CSAV schemes for Case 1. Reference data is taken from \cite{hysing2009quantitative}.}
\label{Case1}
\end{figure}

\begin{figure}[H]
\centering
\subfigure[]{
	\label{case2_Mass}
	\includegraphics[width=73mm]{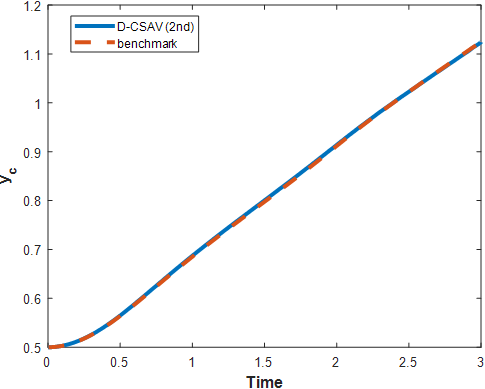}}
\subfigure[]{
	\label{case2_vel}
	\includegraphics[width=75mm]{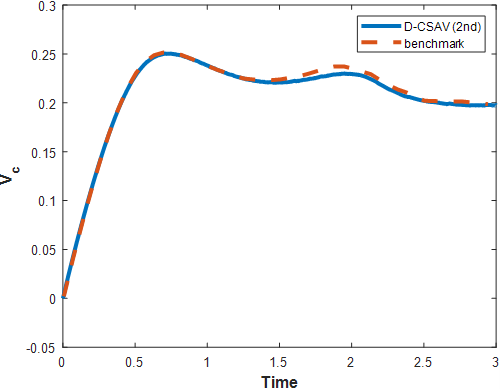}}
\caption{(a) Center of mass and (b) rising velocity obtained using the second order D-CSAV schemes for Case 2. Reference data is taken from \cite{hysing2009quantitative}.}
\label{Case2}
\end{figure}

Figure \ref{Case1} shows the evolution of these benchmark quantities obtained using both the first and second order D-CSAV schemes for Case 1. An excellent agreement with the reference data is achieved for the results of the second order scheme, while a lower accuracy of the first-order scheme is observed. Figure \ref{Case2} presents the evolution of these benchmark quantities obtained using the second order D-CSAV for Case 2. Similarly, the results show good agreement with the reference. { In \cite{wang2024efficient}, they simulated the same problem, where the artificial compressibility scheme was used to discretize the Navier-Stokes equations. A drawback of this method is that it can lead to the velocity oscillation. However, our method does not suffer from this issue.}

\begin{figure}[H]
\centering
\subfigure[$\alpha=10^{-5}$]{
	\label{a_00001}
	\includegraphics[width=75mm]{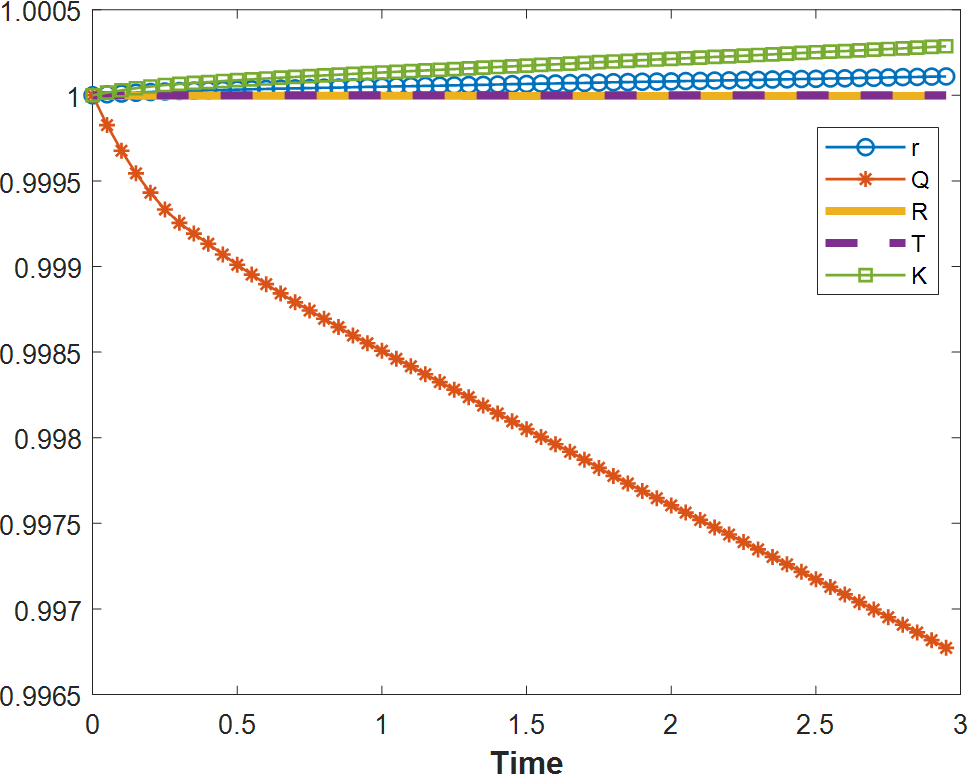}}
\subfigure[$\alpha=1$]{
	\label{a_1}
	\includegraphics[width=73mm]{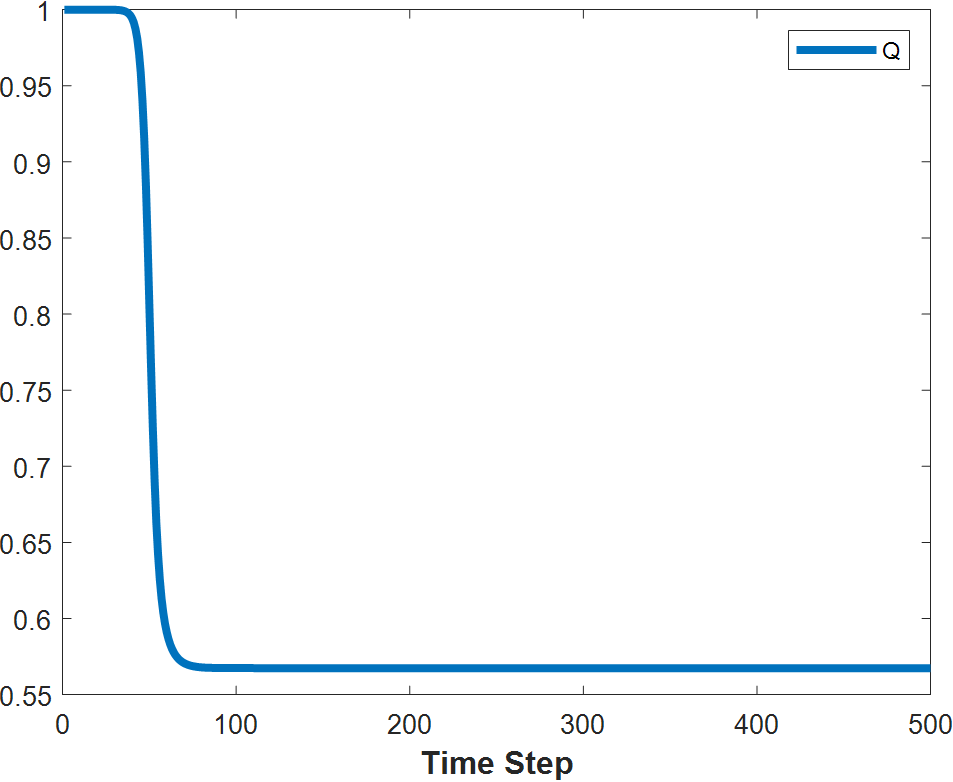}}
\caption{Evolution of scalar auxiliary variables $r$, $Q$, $R$, $K$ using the first order D-CSAV scheme.}
\label{AV_alpha}
\end{figure}

\begin{figure}[H]
\centering
\subfigure[$\alpha=10^{-5}$]{
	\label{a_000012}
	\includegraphics[width=75mm]{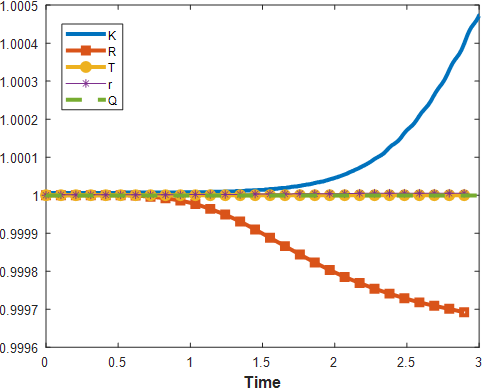}}
\subfigure[$\alpha=1$]{
	\label{a_12}
	\includegraphics[width=73mm]{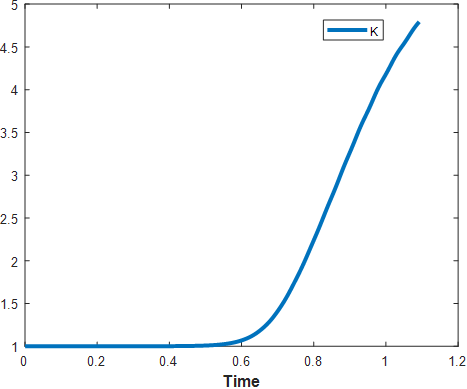}}
\caption{Evolution of scalar auxiliary variables $r$, $Q$, $R$, $K$ using the second order D-CSAV scheme.}
\label{AV_alpha}
\end{figure}

To understand how the parameter $\alpha$ affects the performance of the proposed schemes, we conduct a comparison study using the first order D-CSAV scheme \eqref{CH1_fst}-\eqref{Tn_fst} and the second order D-CSAV scheme \eqref{CH1_sed}-\eqref{Tn_sed} with different values of $\alpha$.
Figure \ref{a_00001} and \ref{a_000012} provide the evolution of the auxiliary variables using the first order scheme and the second order scheme with $\alpha=10^{-5}$. We observe that the values of the auxiliary variables remain close to 1, which confirms the consistency between the modified model and the original model. In Figure \ref{a_1} and \ref{a_12}, we show the evolution of the auxiliary variable $Q$ and $K$ when using $\alpha=1$. The resulting values of these two auxiliary variables $Q$ and $K$ quickly deviate from 1, indicating that the modified model is no longer equivalent to the original model. This comparison result shows that the introduction of $\alpha$ helps maintain the auxiliary variables close to 1; thus, it prevents the potential instability issues of DSAV noted by Remark 2.3 in \cite{fu2021linear}. 
\subsubsection{Bubble merging}
\begin{figure}[H]
\centering
\subfigure[Numerical simulations]{
	\label{3D_coaxial}
	\includegraphics[width=160mm]{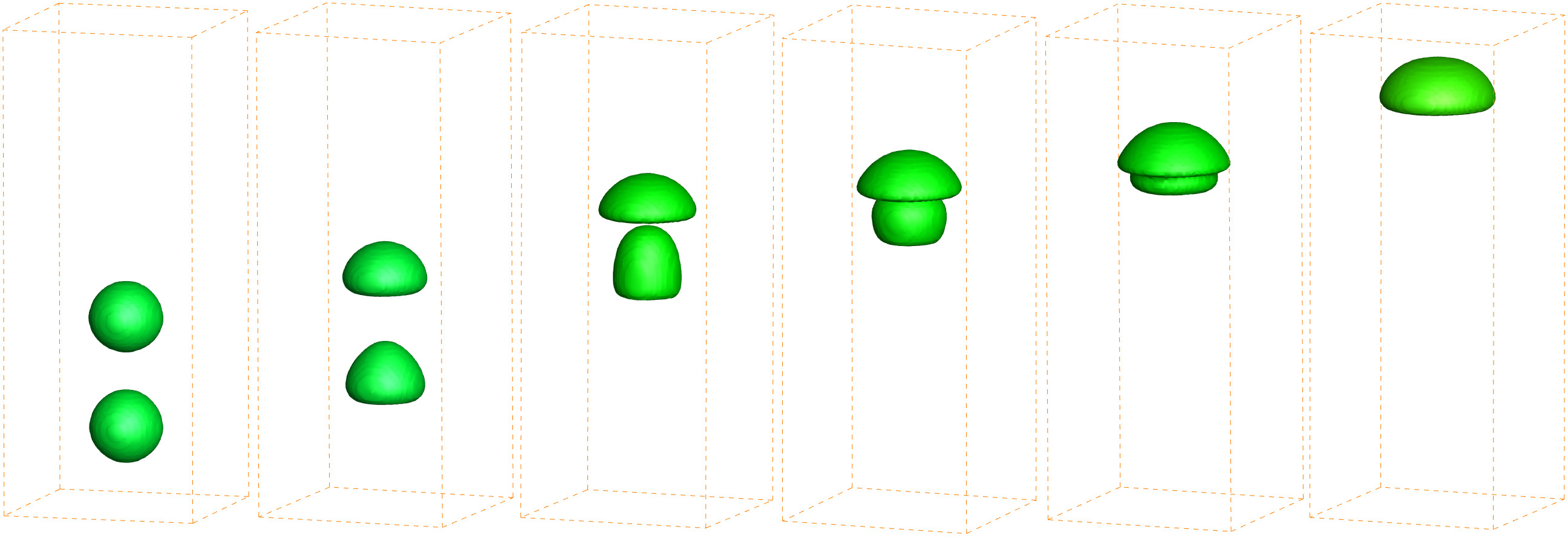}}
\subfigure[Expirimental results]{
	\label{3D_coaxial_exp}
	\includegraphics[width=160mm]{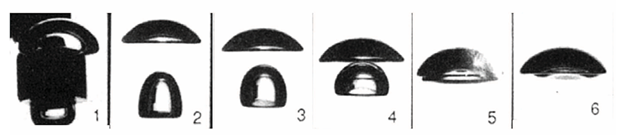}}
\caption{ (a) The dynamical motion of two 3D rising droplets (coaxial coalescence case) under gravity. Profiles of the isosurfaces of the interface { $\phi$ = 0} are taken at $t$ = 0, 0.5, 1, 1.2, 1.5, and 2. (b) The experimental results by Brereton and
 Korotney in \cite{brereton1991coaxial}.}
\label{Case3D_1}
\end{figure}

In this example, we simulate the coalescence dynamics of two 3D droplets under the influence of a gravitational field within a rectangular domain $\Omega=[0,6]\times[0,6]\times[0,12]$. Periodic boundary conditions are applied in the x- and y-directions, while the boundary conditions specified in \eqref{CHNS_bounda} are used for the z-direction.
The model parameters for this case are given as 
\begin{equation}
\begin{aligned}
& \quad g_0=1, \quad \lambda=0.002, \quad \epsilon=0.02, \quad M(\phi)=1, \quad \alpha=0.001, \\
& \rho_1=1, \quad \rho_2=1000, \quad \nu_1=1, \quad \nu_2=100, \quad s=4, \quad \delta t= 10^{-3}.
\end{aligned}
\end{equation}
The initial conditions consist of two droplets of equal size, with one suspended above the other in a coaxial coalescence configuration, described as follows:
\begin{equation}
\begin{aligned}
& \phi^0(x, y, z)=\sum_{i=1}^2 \tanh \left(\frac{0.75-\sqrt{\left(x-x_i\right)^2+\left(y-y_i\right)^2}}{\epsilon}\right)+1, \\
& u^0=(0,0), \quad p^0=0,
\end{aligned} \label{bubble_merge}
\end{equation}
where $x_1=x_2=y_1=y_2=3$, $z_1=1.25$ and $z_2=3.75$. In Figure \ref{3D_coaxial}, we show the isosurfaces of $\phi=0$ at various times. As the two droplets rise under gravity, they exhibit different deformation behaviors. Initially, both droplets take on a cap shape. As they rise, the cross-section of the upper droplet increases while that of the lower droplet decreases. Eventually, they come into contact and merge into one droplet. The numerical simulation shown in Figure \ref{3D_coaxial} qualitatively resembles other numerical results \cite{chen2022highly} and the experiment result in Figure \ref{3D_coaxial_exp}.

\begin{figure}[H]
\centering
\subfigure[Numerical simulations]{
	\label{3D_oblique}
	\includegraphics[width=160mm]{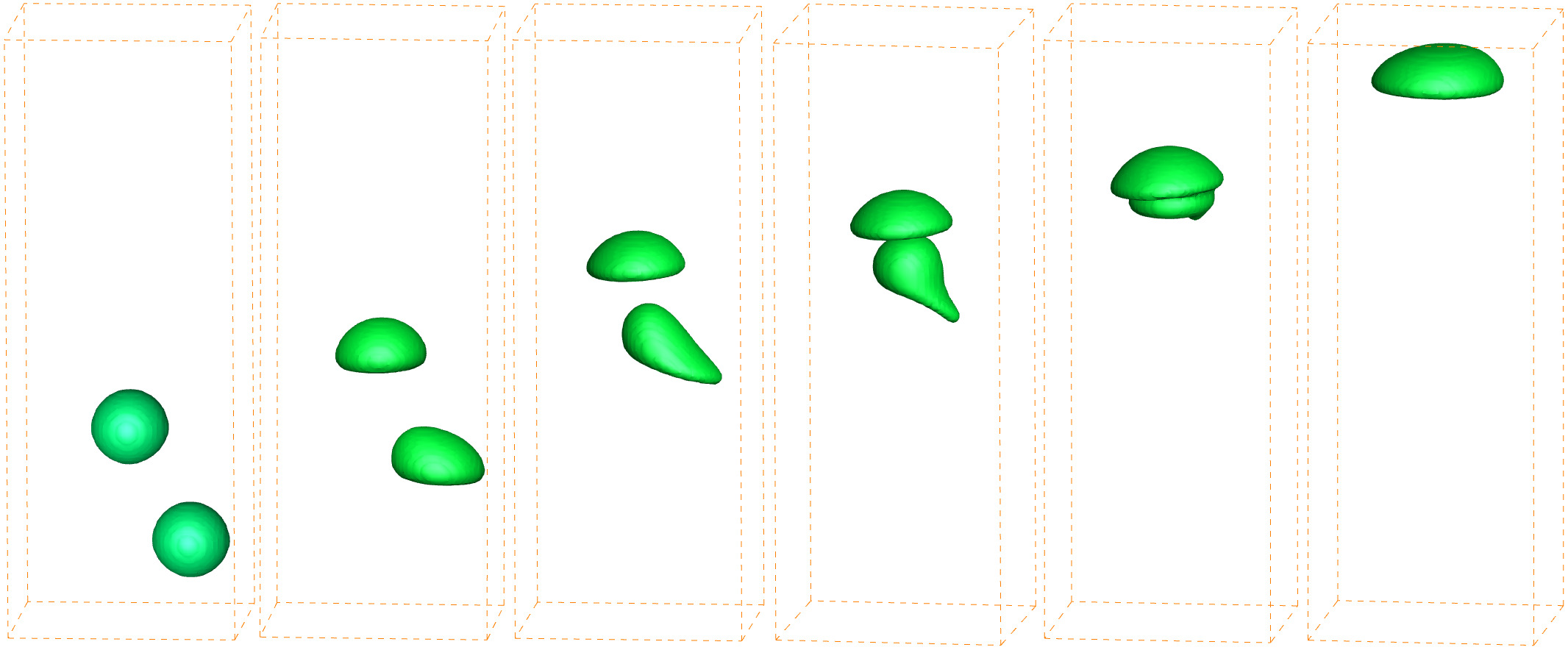}}
\subfigure[Expirimental results]{
	\label{3D_oblique_exp}
	\includegraphics[width=160mm]{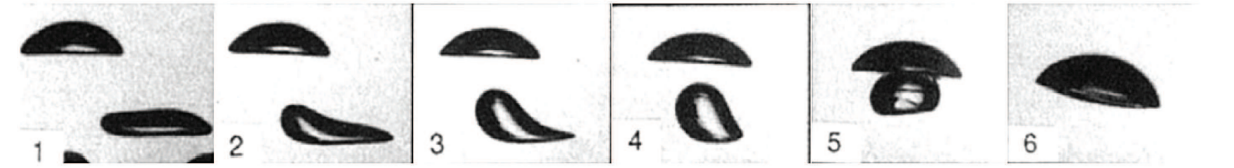}}
\caption{ (a) The dynamical motion of two 3D rising droplets (oblique coalescence case) under gravity. Profiles of the isosurfaces of the interface { $\phi$ = 0} are captured at $t$ = 0, 1.7, 2.3, 2.6, and 3. (b) The experimental results by Brereton and Korotney in \cite{brereton1991coaxial}.}
\label{Case3D_2}
\end{figure}

We further modify the initial conditions to study an oblique coalescence case, where one droplet is positioned staggered relative to the other identical droplet. We reuse the initial conditions specified in \eqref{bubble_merge} but set $\Omega=[0,6]\times[0,6]\times[0,15]$, $x_1=3$ and $x_2=3.6$. The isosurfaces of $\phi=0$ at various times are plotted in Figure \ref{3D_oblique}. Due to the interaction between the two droplets, the upper droplet takes on a slightly inclined cap shape. The lower droplet initially rises towards the upper droplet, then contacts and merges with the upper droplet, creating a larger droplet with a sloping cap shape. We also observe that the numerical results qualitatively agree well with other numerical experiments \cite{chen2022highly} and the experiment results in Figure \ref{3D_oblique_exp}.
\subsubsection{Rayleigh–Taylor instability}
\begin{figure}[H]
	\centering 
	\includegraphics[width=160mm]{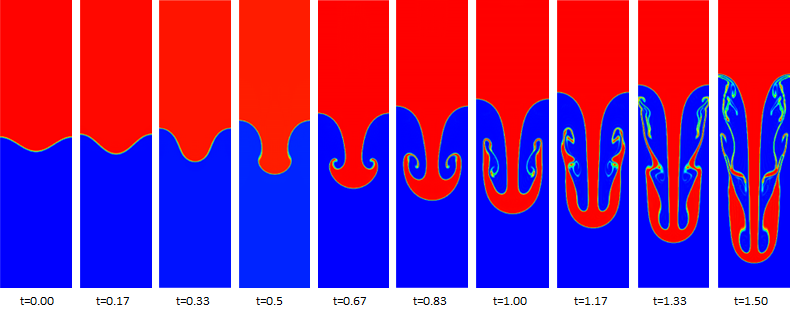} 
	\caption{Evolution of a single wavelength interface to the Rayleigh-Taylor instability.} 
	\label{2D_RT}
\end{figure}

The Rayleigh-Taylor instability occurs in a two-phase system where two fluids of different densities are accelerated relative to each other. In this scenario, if a denser fluid $(F_H)$ is positioned above a lighter fluid $(F_L)$ in a gravitational field, any disturbance at the interface between them can initiate the Rayleigh-Taylor instability. We consider the works of Hosseini \cite{hosseini2017isogeometric} for the validation of our results.
The problem is set up in a 2D rectangular computational domain $[0, d] \times[0,4 d]$. An initial wavy interface separates a lighter fluid in the lower part from a heavier fluid in the upper part of the domain. The initial interface is described by the function
$$
y(x)=2 d+0.1 d \cos (2 \pi x / d)
$$
that represents a planar interface superimposed by a perturbation of wave number $k=1$ and amplitude $0.1 d$. In \cite{hosseini2017isogeometric}, the surface tension coefficient $\sigma$ is set to a small, non-zero value of 0.01. This effectively simplifies the Cahn-Hilliard equation \eqref{CHNS_1}-\eqref{CHNS_2} to a pure transport equation. The other simulation parameters are as follows: $d=1$, $\rho_H=3$, $\rho_L=1$, $\mu_H=\mu_L=0.0031316$ and $g=9.80665$, giving rise to the Reynolds number $Re=\rho_H d^{3 / 2} g^{1 / 2} / \mu_H=3000$. No-slip boundary conditions are applied at the top and bottom boundaries, while free-slip boundary conditions are imposed on the vertical walls. We simulate this problem by using the proposed second-order-in-time numerical scheme \eqref{CH1_sed}-\eqref{Tn_sed}. Figure \ref{2D_RT} illustrates the results for the temporal evolution of the interface over the time interval $[0,1.5]$ with parameters $\delta t=10^{-4}, h=2^{-7}, \varepsilon=0.005$, $M(\phi)=\gamma(\phi^{2}-1)^{2}$, and $\gamma=4\times10^{-5}$.
As anticipated, the heavier fluid on top begins to sink through the lighter fluid, forming spikes that undergo significant deformation. Our numerical results are in good agreement with \cite{hosseini2017isogeometric}.

\begin{figure}[H]
	\centering 
	\includegraphics[width=75mm]{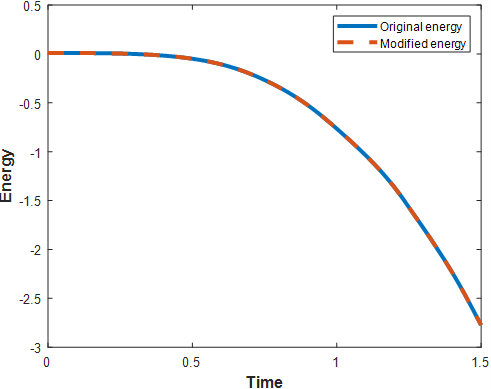} 
	\caption{Evolution of energy in the Rayleigh-Taylor instability simulation.} 
	\label{2D_RT_energy}
\end{figure}

We present the evolution of energy in Figure \ref{2D_RT_energy}, where the impact of gravitational potential energy is taken into account:
\begin{align}
E(\rho, \boldsymbol{u}, \phi)=\int_{\Omega}\left(\frac{1}{2} \rho|\boldsymbol{u}|^2+\lambda \frac{\epsilon}{2}|\nabla \phi|^2+\frac{\lambda}{\epsilon} F(\phi)\right) d \boldsymbol{x}+\int_{0}^{t}\int_{\Omega}\rho\mathbf{g}\cdot\boldsymbol{u}d\boldsymbol{x}d\tau. \label{EO_g}
\end{align}
It is evident that the energy is decaying over time, indicating good energy stablity of the algorithm.

\section{Conclusions}
In this paper, we introduce a second-order-in-time, fully decoupled, and energy-stable numerical scheme for solving the Cahn-Hilliard-Navier-Stokes phase-field model with variable density and viscosity. The scheme is constructed based on the ``zero-energy-contribution'' property and maintaining conservative time discretization for the ``non-zero-energy-contribution'' terms.   By introducing constant scalar auxiliary variables and the associated ODEs,  we have successfully achieved an energy-decaying property. In addition, after numerical discretization, these auxiliary variables are decoupled from other model variables, leading to enhanced computational efficiency. Each time step entails solving three linear elliptic systems, with two of them featuring constant coefficients. We provide a detailed energy stability analysis and conduct extensive simulations involving both 2D and 3D examples. Our results have been compared quantitatively and qualitatively with existing numerical simulations and experimental data, showcasing the efficacy of our proposed scheme.

 \section*{Acknowledgments}
 X.-P. Wang acknowledges support from the National Natural Science Foundation of China (NSFC) (No. 12271461), the key project of NSFC (No. 12131010), Shenzhen Science and Technology Innovation Program (Grant: C10120230046), the Hetao Shenzhen-Hong Kong Science and Technology  Innovation Cooperation Zone Project (No.HZQSWS-KCCYB-2024016) and the University Development Fund from The Chinese University of Hong Kong, Shenzhen (UDF01002028). L. Luo acknowledges support from NSFC 12371442, Macau FDCT 0090/2022/A2, and University of Macau MYRG 2022-00051-FST and SRG 2021-00024-FST.


\bibliographystyle{model1-num-names}
\bibliography{refs}
\biboptions{numbers,sort&compress}


\end{document}